\documentstyle[12pt]{article}
\begin{document}

\newtheorem{lem}{Lemma}[section]
\newtheorem{prop}{Proposition}
\newtheorem{con}{Construction}[section]
\newtheorem{defi}{Definition}[section]
\newcommand{\hf}{\hat{f}}
\newtheorem{fact}{Fact}[section]
\newenvironment{proof}{{\bf Proof:}\newline}{\begin{flushright}$\Box$\end{flushright}}
\newcommand{\ex}{{\cal EX}}
\newcommand{\Cr}{{\bf Cr}}
\newcommand{\dist}{\mbox{dist}}

\title{Quasisymmetric Conjugacies Between Unimodal Maps}
\author{M. Jakobson\thanks{Work supported by the NSF
grant \#o1524481.} \and G. Swiatek\thanks{Supported by NSF grant 431-3604A.}}
\date{September 13, 1991}
\maketitle
\begin{abstract}
It is shown that some topological equivalency classes of
S-unimodal maps are equal to quasisymmetric conjugacy classes.
This includes some infinitely renormalizable polynomials of
unbounded type.
\end{abstract}

\section{Introduction}
\subsection{Quasisymmetric classification of unimodal maps}
\paragraph{Unimodal maps.}
We discuss unimodal maps of the interval. A standard example is the
quadratic family 
\[ x\rightarrow ax(1-x) \]   
where $a$ is a parameter from the interval $(0,4]$. 
Other important classes are maps extendable in an analytic
quadratic-like fashion in the sense of \cite{dohu} and the S-unimodal
class where no analytic extension is postulated, instead the map is
assumed to have negative Schwarzian derivative. In the following
discussion, unless otherwise indicated, we mean maps from the union
of these two classes.
  
Unimodal maps exhibit impressively rich dynamics. The framework for
studying them was laid by \cite{mithu}. There, topological dynamics of
unimodal was described in terms of the {\em kneading sequence}. However,
the basic idea the kneading invariant can be traced back to an earlier
paper \cite{metro}.  

In some cases, the dynamics of unimodal maps has been well understood.
This includes maps with periodic or preperiodic kneading sequences for
which the analytic cases were studied in \cite{mss} and \cite{etiuda}.
In this paper, we confine ourselves to other, or aperiodic, invariants.   

\paragraph{Quasisymmetric conjugacies.}
By the work of \cite{guke1} we know that two maps with the same
aperiodic kneading sequence are topologically conjugate. The
conjugating homeomorphism is quasisymmetric if and only if it can be
extended to a quasiconformal hoemomorphism of the plane. 
By the celebrated theorem of \cite{abu}, this is equivalent to the ratio 
\[ \frac{g(x+h)-g(x)}{g(x)-g(x-h)} \] 
being uniformly bounded for all real $x$ and $h$ so that the
relevant points are in the domain of $g$. Moreover, there exists a
quasiconformal extension whose norm is bounded in a uniform way in
terms of the supremum of this ratio (which we will call the
quasisymmetric norm.)

It is known that quasisymmetric homeomorphisms are H\"{o}lder
continuous, but usually not absolutely continuous. 

In dynamics, the idea of studying quasiconformal (-symmetric)
conjugacy classes has been introduced and proven stunningly successful
by a series of works by D. Sullivan. A recent work \cite{miszczu}
deals directly with unimodal maps and was an inspiration, as well as
the starting point of this work.   

\paragraph{Various results on quasisymmetric classification.}
The standing conjecture is that the quasisymmetric conjugacy classes
are equal to topological conjugacy classes for aperiodic kneading
invariants. This conjecture has so far been proven in three cases.

First, there is a class of infinitely renormalizable maps which was
treated in \cite{miszczu}. 

Secondly, in the Misiurewicz case which means that the critical point is
not recurrent the conjecture was proved by \cite{mike}.

Finally, a recent result of Yoccoz should be mentioned
which implies the conjecture for all non-renormalizable polynomials in
the analytic polynomial-like class. This work has not yet been
circulated; the reader may, however, consult \cite{hub}. 

\paragraph{What this paper contributes.}
We prove the conjecture for some, not all, non-renormalizable maps in
the S-unimodal class. We also show a new approach to renormalizable
cases. We prove the conjecture in some infinitely renormalizable cases
where it is new even in the polynomial class.

\paragraph{Consequences of our results.}
A famous consequence is that 
is that if the conjecture is proven for any kneading sequence in the
polynomial class, the corresponding component of the Mandelbrot set
reduces to a point. This is proven by the pull-back construction of 
\cite{miszczu} and a deformation argument of the kind used in
\cite{mss} and \cite{sula}. 

Another consequence concerns the existence of absolutely continuous invariant
measures. The Collet-Eckmann condition (see \cite{cole}) is shown to
be a topological invariant of non-renormalizable maps in our class.
The proof of this remark is in Section~\ref{sec:5,3}. 

The question of topological invariance of the Collet-Eckmann condition
in the class of S-unimodal non-renormalizable maps was stated by 
J. Guckenheimer \cite{conf}. There is a related question of whether
the existence of an absolutely continuous invariant measure is a
topological property for S-unimodal maps. This, however, seems harder
and at this point we can only state it as a problem.   

\paragraph{Acknowledgments.}
MJ acknowledges hospitality of the Thomas B. Watson IBM Center and
Institut des Hautes Etudes Scientifiques where parts of this work were done.

\subsection{Induced maps}
\paragraph{Assumptions.}
The class of functions ${\cal C}$ is defined by the conditions:
\begin{defi}\label{defi:n5,1}
\begin{enumerate}
\item
Each $f\in {\cal C}$ maps the interval  $[-1,1]$ into itself.
\item \label{as 2}
Functions from $\cal C$ are three times differentiable and, wherever the 
first derivative is nonzero, their Schwarzian derivative is non-positive.
\item
Each function $f\in {\cal C}$ can be represented as 
$h(x^{2})$ with $h$ being a diffeomorphism.
\item \label {as 4}
The critical value $h(0)$ is greater than $0$.
\end{enumerate}
\end{defi}

These assumptions in particular imply that $f(-a)=f(a)$.

\begin{defi}\label{defi:n5,4}
Let $\alpha$ be an aperiodic kneading sequence. Then ${\cal C}_{\alpha}$ is
defined to be the set of all maps from $\cal C$ with this kneading sequence.
\end{defi}

In this paper, we only deal with maps whose kneading sequence is
aperiodic.

Assumption ~\ref{as 4} implies that there exists a fixed point $q$ for every $f\in {\cal C}$ with $q>0$. 
We consider the induced map $\phi$ defined to be the first return map
on the interval between $q$ and $-q$. This interval will be called the
{\em fundamental inducing domain} of $f$.

It is easy to see that
the induced map consists of a number of continuous branches all of which except
one are monotonic. Also, the construction is topological, by which we mean that
if maps $f_{1}$ and $f_{2}$ are topologically conjugate, the same is true of 
their induced maps.

\begin{defi}
Given an interval $I\subset [-1,1]$ we define a {\bf stopping rule} on 
$I$ to be a continuous positive integer valued function defined on 
an open subset of $I$.
\end{defi}

\begin{defi}\label{defi2,1}
An {\bf induced map} of $f\in {\cal C}$  on an interval $I\in [-1,1]$ is a map 
of the form 
\[x\rightarrow f^{s(x)}(x)\]
where $s(x)$ is a stopping rule on $I$ and we mean that the induced map is not
defined where the stopping rule is not.
\end{defi}

So, induced maps and stopping rules are really the same thing and we will
keep in mind that one always determines the other. 

\begin{defi}\label{defi:9,2}
A restriction of an induced map $\varphi$ to a connected component of its 
domain will be called a {\bf branch} of  $\varphi$.
\end{defi}

\begin{defi}
An {\bf induced monotone branch} is an induced map with a constant stopping 
rule whose domain is an interval and which is monotone.
\end{defi}

\begin{defi}\label{defi:10,1}
An induced monotone branch defined on an interval $(a,b)$ is said to be
$\epsilon$-{\bf extendable}
if there is an induced monotone branch $g$ with the same stopping rule defined
on a larger interval $(c,d)\supset (a,b)$ such that the cross-ratio
\[\frac{|g(a)-g(c)||g(b)-g(d)|}{|g(a)-g(d)||g(b)-g(c)|}\]
is more than $\epsilon$.
\end{defi}

In the future, we will fix a uniform value of $\epsilon$ and simply talk of 
extendable maps. Monotone extendable branches have bounded distortion
(see \cite{guke2}).

\begin{defi}\label{defi:9,1}
A {\bf critical branch} is a branch of the form $g(x^{2})$ where $g$
is a monotone branch, defined on a symmetric neighborhood of $0$.  
\end{defi}

Here, it is understood that the domain of $g$ may very well be larger
than the image of the domain of the map by the quadratic map. Hence,
our notion of the {\em image} of a critical branch is non-standard, as we
define it to be the image of $g$.  

\begin{defi}\label{defi:9,3}
A critical branch $g(x^{2})$ is {\bf extendable} if $g$ is. 
\end{defi}

\begin{defi}\label{defi:9,4}
A branch with domain $P$ and stopping rule $s$ is said to be {\bf folding} 
(extendable) if there is an $\overline{s}<s$ such that 
\begin{itemize}
\item
$f^{\overline{s}}$ on $P$ is an induced monotone branch (extendable).
\item
$f^{s-\overline{s}}$ on $f^{\overline{s}}(P)$ is a critical branch 
(extendable).
\end{itemize}
\end{defi}

The image of a folding branch is, by definition, equal to the image of
the corresponding critical branches.

\begin{defi}\label{defi:q21a,1}
If $\xi$ is a branch with stopping time $s$, a {\em settled branch} can be
defined for any settling time $\overline{s}\leq s$. The settled branch
is always $f^{\overline{s}}$, and its domain is equal to the domain of
$\xi$. If $f^{\overline{s}}$ folds on the domain of $\xi$, we also
need to specify the image. By definition, it is equal to the preimage of the
image of $\xi$ by $f^{s-\overline{s}}$ (which is well-defined since
$f^{s-\overline{s}}$ is invertible on the relevant interval).   
\end{defi}

\paragraph{Prefered induced maps.}
We describe a class of induced maps which have particularly useful properties.

\begin{defi}\label{defi 1,20}
An induced map is called a {\bf preferred map} if it has the following
properties:
\begin{itemize}
\item
All branches are either monotone or folding and extendable.
\item
All folding branches have the same critical value whose image is not entirely
contained in one the external branches.
\item
The branches do not accumulate at the endpoints of the interval, and the 
external branches are monotone.
\end{itemize}
\end{defi}

\subparagraph{Notational conventions.}
We will use parallelism in our notations between objects defined for $f$ and
similarly defined objects for $\hf$ which automatically receive the same
labeling only marked with a hat $\hat{ }$ sign.  

Another problem comes from a considerable number of uniform constants which
will abound in future arguments. To say that a constant is ``uniform''
means that it depends only on global distortion properties of maps $f$
and $\hf$. More 
precisely, it only depends on the $C^{3}$ norm of the corresponding map $h$
and the infimum of derivative of $h$.
Uniform constants will be denoted by the letter $K$ with a subscript.

A statement which contains uniform constants should mean that ``for
each occurrence of a uniform constant, there exists a uniform numerical
value which makes the statement true.'' We do not claim that uniform
constants denoted with the same letter correspond to a fixed value
throughout the paper. Thus, $K_{1} > K_{1}$ would be considered a true
statement, though we will use subscripts to avoid such extreme
examples.

\subsection{Non-renormalizable maps of basic type}

\paragraph{Two maps.}
From now on we consider a pair of maps, $f$ and $\hf$, both from $\cal C$ with
the same kneading sequence. It is known that under our assumptions they are
topologically conjugated so that
\[\hf=h^{-1}\circ f\circ h\; .\]

In this context, we can talk about {\bf equivalent stopping rules} $s$ and 
$\hat{s}$ if the relation if the domain of $s$ is mapped onto the
domain of $\hat{s}$ by the topological conjugacy $H$, and 
\[s=\hat{s}\circ H\]
holds where defined.

\paragraph{Basic construction.}

The way we refer to the topological dynamics of our maps is through the
basic construction as it stands in \cite{yours}.

We assume that
 
\begin{em}
On each stage of the construction, the critical value falls into a monotone
branch.
\end{em}  

Since the basic construction is topological, this is a topological
condition. We will refer to it as ``basic dynamics.'' 

We do not know how to express this assumption in the language of
kneading sequences. However, we note that the basic class is wider
that the set discussed in \cite{jak1}. On the other hand, in the
complex analytic case it is narrower than the
intermittently recurrent class considered by Yoccoz in his recent
work. 
 
\paragraph{The main result.}

{\bf Theorem 1}

\begin{em} 
Any two maps from $\cal C$ with basic dynamics are quasisymmetrically 
conjugate;
moreover, the quasisymmetric norm of
conjugacy is bounded by a uniform constant.
\end{em}

\paragraph{Two important structures of a folding map.}

Any map from $\cal C$ has two distinguished points: the critical value
and the fixed point $q$ inside $(-1,1)$. Importance of the forward
critical orbit is well-known. In particular, its combinatorics defines
the topological class of the map. 

However, there is another structure worth looking at, and that is the
backward orbit of the fixed point. In non-renormalizable cases this
orbit is dense, thus any homeomorphism which maps backward preimages
of the fixed point of one map onto corresponding points of another map
with the same dynamics must be the conjugacy. This is how we build the
conjugacy in this work as a limit point of ``branchwise equivalences''. 

The forward critical orbit continues to play an important role in our
construction, and the reader may note how both concepts interact in
our ``critical pull-back'' and ``marking'' operations.  
    
\paragraph{Introducing branchwise equivalences.}

\begin{defi}\label{defi:q10a,1}
A {\bf branchwise equivalence} is a triple which comprises two equivalent
stopping rules together with a homeomorphism which maps the domains of branches
of one map onto the corresponding domains of branches of the other
map.
\end{defi}

The homeomorphism from the domain of $f$ to the domain of $\hf$ which
is the third component of the branchwise equivalence will also be called
a branchwise equivalence, and the induced maps will then be referred
to as the ``underlying induced maps''.  

The subset of $I$ on which a branchwise equivalence coincides with the
conjugacy will be called its {\em marked set}. By definition, the
marked set contains at least the boundary of the domain of $s$.

The domains of branches of $f^{s}$ as well as components of the
interior of the complement of the domain of $s$ will be called the
{\em domains} of this branchwise equivalence. We will thus speak of
monotone and folding domains, while the last kind will be referred to
as {\em indifferent domains}.  

A branchwise equivalence which has no indifferent domains (i.e., the
domain of the underlying stopping rule $s$ is dense) will be called
{\em regular}.

The construction that we use up to Section~\ref{sec:6} only gives
regular branchwise equivalences. 
\paragraph{The strategy of the proof.}

With every inducing construction there is an associated procedure of refining
the domains of branches. This gives a natural way of building up a conjugacy
between two maps. With any pair of equivalent induced maps we can 
associate a branchwise equivalence. If our inducing construction is 
sufficiently general, as it is the case with the basic construction followed
by inducing on all monotone branches, we may
hope that the actual topological conjugacy can be found somewhere in the 
closure of these branchwise equivalences. If so, the only thing remaining 
is to show that all branchwise equivalences in the class we consider are
uniformly quasisymmetric. This is, of course, the hardest part.

Our basic technique will be patching different branchwise equivalences
together to get new branchwise equivalences. This kind of procedure cannot
be effectively carried out using real maps only. This is one reason why we
will complexify our problem and indeed work on the level of quasiconformal
extensions of branchwise equivalences. 

So first, we are going to define the procedure of inducing and at the same
time of constructing branchwise equivalences. We will then check that that 
construction is sufficiently general, so that in fact the basic construction
can be approached using our methods.

Then, real work begins. We will redefine the construction in terms of 
quasiconformal extensions of branchwise equivalences. The complex 
procedure will be designed so as to guarantee that complex 
quasiconformal norms of the maps we construct will be uniformly bounded.  

\subsection{Renormalizable maps}
\paragraph{Statement.}
Let $f$ now be a renormalizable map. Let $I^{i}$, $i>0$, be the maximal
decreasing sequence of restrictive intervals around $0$. Then
$I_{j}^{i}$ denote the orbit of $I^{i}$ by the map.  

We also get a sequence of maps from $\cal C$. Here, $f_{0}:=f$ and
$f_{i}$ is the first return map on $I^{i}$ conjugated by an affine map
so that $I^{i}$ becomes $[-1,1]$. 

We will say that a preferred regular induced map of a renormalizable
map is  suitable if
one of $I^{1}_{j}$ is in the domain of a folding branch and mapped
into itself by the branch.   

Our main result here is that:

{\bf Theorem 2}

\begin{em}
Consider two infinitely renormalizable topologically conjugated maps
from $\cal C$. Suppose that for each $i$ there exist a quasisymmetric
branchwise equivalence between suitable preferred induced maps of
$f_{i}$ and $\hf_{i}$. 
If their quasisymmetric norms are uniformly bounded, then $f$ and
$\hf$ are quasisymmetrically conjugate and the qs norm of the
conjugacy is bounded by a uniform function of the common bound.   
\end{em}

\paragraph{Comment.}
Theorem 2 may be applied in various situations. If the infinitely
renormalizable map has ``bounded type'' as introduced in
\cite{miszczu}, the assumption is relatively easy to verify, because
all suitable maps are finitely complicated and subject to ``bounded
geometry''. 

However, we hope that usefulness of Theorem 2 extends far beyond that.
The reader recalls that our approach to the conjugacy in
non-renormalizable cases is by building more and more refined
branchwise equivalences. The point is that quite often the suitable 
equivalence can be constructed in our way, and really the fact that
the map is renormalizable makes no difference in the construction.
Hence, a uniform bound on the quasisymmetric norms follows
exactly as in the renormalizable case. 

That means, for example, that if all suitable induced maps can be
obtained in the basic construction, the conjugacy is uniformly
quasisymmetric. This defines a class of infinitely renormalizable maps
for which, as far as we are aware, the quasisymmetric conjugacy
result is new even in the polynomial case. 

Moreover, from the point of view of our approach of refined branchwise
equivalence, all that is needed to close the infinitely renormalizable
case is a uniform estimate for the construction in remaining
non-renormalizable cases, called ``box cases'' in \cite{yours}.     

Section~\ref{sec:6} also contains a theorem for finitely
renormalizable maps, which we will not discuss here.

\section{Branchwise equivalences.}
\subsection{Introduction of branchwise equivalences.}
We assume that we are given two maps $f$ and $\hat{f}$ which satisfy our
assumptions. 

We will give the description of our construction as a recursive procedure.
That is, we are going to show simple primary objects and define operations
allowing us to construct more complicated things from those simplest ones.
From now on, we assume that the basic dynamical intervals $(-q,q)$ and
$(-\hat{q},\hat{q})$ have been uniformized by the affine maps from the unit
interval. So we will simply assume they are both $[0,1]$.

Before we continue, we would like to make a comment on the condition that
images of folding branches must not be contained in an external branch.
A simple observation is that if all other conditions for the map to be preferred
are satisfied except for this one, there is a simple way ``adjust'' the map to
become preferred. Namely, we can compose folding branches with that external 
branch. Since a repelling (pre-)periodic point is an  endpoint of the external 
branch, the image of the fold will become longer. If we continue to compose
until the critical value leaves the domain of the external branch, we will
get a preferred map. 
\subparagraph{Boundary-refinement.}
There is a typical construction which we now describe. We can consider the
leftmost branch of the map and compose it with the map itself. Then, we can
take the new leftmost branch and again compose it with the original map.
As this procedure is repeated, the leftmost branch gets exponentially shorter.
If we continue the process to infinity, we get something that will be called
{\em a map infinitely boundary-refined on the left.} Of course, we can also
construct maps infinitely boundary-refined refined on the right or on both 
sides. We could also choose a point $x$ very close to $0$ and continue the
left boundary refinement until $x$ is no longer in the leftmost branch. This 
would be the {\em boundary refinement to the depth of $x$}. If we start with
equivalent induced maps, and the depth of the refinement is determined by 
topologically equivalent points, then the resulting maps will also be 
equivalent.

\subparagraph{Boundary refinement of a branch.}
Suppose we are given a preferred induced map. Any {\em monotone} branch of 
this map can be boundary-refined as follows: we first boundary-refine the
whole map, and then compose the branch with the result of the refinement. 
Everything we said of refining a map has obvious consequences for this 
construction. However, there is one particular case we want to discuss.
Suppose our monotone branch shares its right endpoint with a very short 
folding branch. Very short critical means that the critical value falls into 
one of the external branches and takes a long time to leave the external 
branch adjacent to $1$.\footnote{$1$ is $q$ after reparametrization: a 
repelling fixed point for the induced map.} 

We will often want to refine the monotone branch so that the domain of the 
rightmost branch of the result has length uniformly comparable with the 
length of the adjacent folding branch. If the critical value takes $n$ iterates
to leave to rightmost branch, we choose a point in the domain of the monotone 
branch which hits also stays in the domain of the rightmost branch for $n$ 
iterates, and its $n$-th image hits the boundary point of that branch. 
If we then refine the monotone branch to the depth of the image of this point,
then, indeed, the length of the domain of the new branch adjacent to the 
folding branch is comparable with the length of the domain of the folding 
branch itself. Also, this construction is topological: if we start with a pair
of equivalent maps, we get equivalent maps.

We call it the {\em refinement to the depth of the adjacent folding branch.}

 \subsection{How to build branchwise equivalences.}
\paragraph{The primary branchwise equivalence.}
We need a preferred branchwise equivalence to begin with. We would also like
it to satisfy two estimates uniform with respect to the choice of $f$ and 
$\hf$.
\begin{enumerate}
\item
For either map, the lengths of the domains of any two adjacent branches are
comparable, i.e. their ratio is bounded and bounded away from $0$ by a uniform
constant. 
\item
The primary branchwise equivalence is affine inside the domain
of any branch, and the identity outside the interval $[-1,1]$.
\item
The quasisymmetric  norm of the primary branchwise equivalence is
uniformly bounded.
\end{enumerate}

\paragraph{Finding the primary branchwise equivalence.} 
A reasonable way to start is by looking at the induced maps $\phi$ and 
$\hat{\phi}$ as defined in the introduction. They are always preferred. 
However, the additional estimates may not hold and the only way that can 
happen is when the central folding branch is
extremely short and, as a result, the two monotone branches adjacent to it
are non-extendable. Each of these non-extendable branches maps onto the whole 
interval $[0,1]$, and, thus, they can be composed with the original map,
or ``refined'' as we will often refer to this procedure. 

First, we consider the situation when there are at least two monotone branches
on each side of the central folding branch. 
Since the situation is wholly symmetric, it is enough to analyze what happens 
to the non-extendable branch on the left of the central folding branch.
We refine it on the right to the depth of the central critical branch.
Now, the branch adjacent to the central branch will already be extendable.
But the refinement procedure will also create preimages of the central
branch together with its non-extendable neighbors. 

So in the next step, we start we our original map again and this time refine
the non-extendable branches by composing them with a suitably boundary-refined
version of the map obtained on the previous step. As we continue the process to
infinity, the non-extendable branches eventually are crammed into a Cantor 
set (of zero measure.) Note also, that the procedure leaves the external 
branches of the original map unaffected. 

However, if initially there is only one monotone branch on each side, this
procedure would lead to a non-preferred induced map in which branches 
accumulate to $0$ and $1$.

Fortunately, there is another solution available in this case. If the central
branch is very short, the induced map intersects the diagonal at a point
$q'$, which is repelling with period $2$. We now consider 
the first return map from the interval $-q',q'$ onto itself. This turns out
to a preferred-map with one folding and infinitely many monotone branches.
Hence, the previously described construction applies.  

\paragraph{Various types of primary branchwise equivalences.}
In the previous paragraph we indicated how to construct a pair of initial
preferred induced maps. Once we have them, we can construct a branchwise 
equivalence between them. There will be two kinds of branchwise equivalences:
marked and unmarked. The unmarked branchwise equivalence by the definition is
the identity outside of the interval $[0,1]$ and is affine in domains of all
the branches. The marked equivalence is so designed as to map a point in the
domain of one of the monotone branches to its corresponding point under
the topological conjugacy. Typical the point here will belong to the forward 
critical orbit. To achieve that, we make the branchwise equivalence
linear fractional (for example) on the domain of this branch.

\subparagraph{Boundary-refined versions.}
The map defined above can also occur in infinitely many boundary-refined 
versions. The main technical problem that we encounter with boundary-
refined branchwise equivalences is that the banal extension by the identity
beyond the unit interval does not work. Simply, as we consider boundary-
refinements to growing depth, the quasisymmetric norm deteriorates (unless both
$f$ and $f'$ have the same eigenvalue in $q$), and for the infinite depth 
boundary-refinement this extension could not be quasisymmetric at all.

So we use another extension. Let us say that we want to obtain the right
boundary-refinement of the primary map. We do the boundary-refinement as
previously described and construct the branchwise equivalence in the same way
we showed in the last paragraph. Then, we extend by the identity to the left
of the unit interval and mirror the result about $1$ to extend it to the right
of $1$.

This gives us a quasisymmetric map. 

\subparagraph{Summary of primary branchwise equivalences.}

Our future estimates will depend on the following estimates for this primary 
map:
\begin{itemize}
\item
The maximum ratio of lengths of any two adjacent domains.
\item
The ratio of $1$ to the lengths of the external branches in case of branchwise
equivalences which are not boundary-refined.
\item
The maximum quasisymmetric norm of any branchwise equivalence.
\end{itemize}

\subsection{Branchwise equivalences build-up}

Our next task is to describe how to ``refine'' primary branchwise equivalences
so as to make them approach the actual conjugacy. The process is more or less
parallel to the basic construction. However, one important difference is that
we want uniformly quasisymmetric maps on all stages of the construction, which
something that the basic construction does not provide.

We first describe how to obtain branchwise equivalences which are not boundary-
refined.

\paragraph{The first operation: critical pull-back.}
To perform this operation we need a branchwise equivalence 
$\Upsilon_{1}$. For the underlying pair of induced maps, we pick a pair
of corresponding folding branches $\psi$ and $\hat{\psi}$. We assume that these
branches are extendable, and that their images are not contained in an external
branch. We also need another branchwise  equivalence $\Upsilon_{2}$ such that 
the critical value
of $\psi$ falls into a monotone branch of $\Upsilon_{2}$.
\footnote{We will often use the word ``branch'' to really mean ``the domain
of a branch''. We hope that it will not lead to confusion, while making our
text smoother. } 

First, we will describe what happens on the level of the associated 
stopping rules. In terms of induced maps critical pull-back can be expressed
as composing $\psi$ with the induced map associated with $\Upsilon_{2}$ and 
the same thing with $\hat{\psi}$ in the other map. However,
here is one exception: if the critical value of $\psi$ is in one of the 
extreme domains of $\Upsilon_{2}$, then we continue composing the resulting 
critical branch with this extreme branch until the critical value leaves 
its domain. This, however, will not create any new branches.

Then we proceed to define the construction of the branchwise equivalence between
the new induced maps.

First, we take a $\em marked$ version of $\Upsilon_{2}$. Namely, we require
that  
\[c(\hat{\psi})=\Upsilon_{2}\circ c(\psi)\]
where the notation $c(\cdot)$ means ``the critical value of''.   
Marking means changing the branchwise equivalence $\Upsilon_{2}$. We
will later describe this process precisely. Right now we only assume
that marking does not alter $\Upsilon_{2}$ except on the monotone
domain  which contains the critical value.

The marking ensures that $\Upsilon_{2}$ can be lifted by $\psi$
and $\hat{\psi}$ and, by the definition, it is the order-preserving lift  
that is going to replace $\Upsilon_{1}$ inside the domain of $\psi$.

Outside the domain of $\psi$, the branchwise equivalence is left unchanged. 

Please note that if we have a number of folding branches, the critical 
refinement on one of them commutes with the critical refinement on any 
other, and in this sense we can say that the critical refinement can be
done concurrently on all folding branches of a given map.
\paragraph{The second operation: monotone pull-back.}

We need a branchwise 
equivalence $\Upsilon_{1}$ and a pair of corresponding extendable monotone 
branches $\xi$ and $\hat{\xi}$ of the underlying induced maps. We also need 
another branchwise equivalence $\Upsilon_{2}$. 

On the level of induced maps, this operation is simply composing $\xi$
and $\hat{\xi}$ with corresponding maps associated with $\Upsilon_{2}$.

To get the new branchwise equivalence, we replace $\Upsilon_{1}$ on the domain
 of $\xi$ with
$\hat{\xi}^{-1}\circ\Upsilon_{2}\circ\xi$ and leave it alone outside of the
domain of $\xi$.

\paragraph{Boundary refined versions.}

If we use $\Upsilon_{1}$ which is not boundary-refined in a pull-back step, 
then we get a non-boundary-refined map as a result. To obtain a 
boundary-refined version of the 
result to a certain depth, we should start with $\Upsilon_{1}$ refined to this
depth.

\paragraph{A step of the construction.}
In the preceding paragraph, we described basic elements of the construction.
Now, we will show how a step of the basic construction can be mimicked using
these techniques.

Our starting point is a preferred branchwise equivalence $\Upsilon$ (we also
know how to construct its marked versions.) Our construction will eventually
yield another preferred branchwise equivalence, and the change of the folding
branches will change in the same way as in the basic construction. However,
we will watch that these properties are satisfied, which are considered 
unimportant in the basic construction:
\begin{itemize} 
\item
The lengths of any two adjacent domains are comparable.
\item
If two monotone branches share an endpoint, they are always refined 
simultaneously, and if fact they are subject to an infinite-depth boundary
refinement at their common endpoint.
\end{itemize}

\subparagraph{The first critical pull-back and boundary refinement.}
By assumption, the critical value is in a monotone branch. If it is
not too  close to an endpoint of the domain, or the adjacent branch is
folding, we simply apply critical pull-back
on all folding branches. If, however, the critical value is too close to an
endpoint, that gives us a map with non-extendable branches. If these non-
extendable branches are monotone, we use another procedure. We boundary-refine
the branch adjacent to the branch containing the critical value. The boundary
refinement has the depth comparable to the distance of the critical value
from the endpoint on one side. What happens on the other side depends on 
whether the next branch is folding. If it is, there is no additional refinement
on that side and this step is completed. If it is monotone, there an infinite-
depth boundary refinement on that side. In fact, all consecutive monotone 
branches are subject to the infinite-depth boundary refinement which will only 
end when a folding branch is encountered. Then after this process has been
completed, the resulting map is pulled-back on all folding branches of
$\Upsilon$, just like in easier case.
 
\subparagraph{The filling-in.}
The previous step resulted in a non-preferred map, because the folding branches
may have different critical values, and also because some folding branches may
not be extendable. The following procedure of ``filling-in'' is completely 
analogous to what was described under the same heading in the basic 
construction. If the result of the first critical pull-back is denoted with
$\Upsilon'$, the second step differs from the first in the the map we pull-back
in $\Upsilon'$, not $\Upsilon$. If a preliminary sequence of boundary-
refinements is needed, we do it on $\Upsilon'$ just like we did on $\Upsilon$ in
the first step. Then, the resulting $\Upsilon''$ is again pulled-back on the
folding branches of $\Upsilon$ and so on to infinity. Eventually, all folding
branches with the critical value in the old place will disappear, squeezed into
a Cantor set, and we get a preferred branchwise equivalence again.

\paragraph{The final refinement of the monotone branches.}
Our objective is to get a sequence of branchwise equivalences which tend to
the topological conjugacy. If only do basic steps as described above, this will
not be the case, since some monotone branches, for example the external ones,
will never be refined. That is why at some moment we have to stop and refine
these ``lagging'' monotone branches. 

This is a little bit similar to the boundary-refinement sequence described 
in the previous paragraph. We refine all monotone branches which are not 
contained in the external primary branches. We all also refine some branches
contained in the external primary branches if it is necessary to preserve
the rules of the chain boundary refinement. At the points of
tangency with folding branches, we refine to the depth so chosen
that the new adjacent branch is always shorter than the folding branch, but 
still in a uniformly bounded away from $0$ ratio. Between adjacent monotone
branches the refinement is infinitely deep. In the next step, we refine all new
monotone branches except for those on which the previous refinement stopped
(i.e. the ones adjacent to the folding branches of the original map.) As we
repeat this process, we eventually wind up with monotone branches all smaller 
than the longest folding branch of the original map. But, that tends to zero
in the basic construction, so indeed we get a sequence of maps that tend to 
the conjugacy on the set between the primary external branches.

\paragraph{Marking conditions.}
\begin{defi}\label{defi:n6,6}
A {\em marking condition} is a choice of an infinite ray which starts in
a monotone domain of $f$.
\end{defi}

\subparagraph{Reduced versions of branchwise equivalences.} 
Given a ray and a branchwise equivalence $\upsilon$ constructed in the
way  just 
described, we want to consider a reduced version of this equivalence
with respect to the ray. To define the ``reduced version'' we need to
exactly mean a
branchwise equivalence which coincides with $\upsilon$ on all branches
of some primary branchwise equivalence which are intersected by the
ray, why all other branches of that branchwise equivalence have been
refined at most twice.

It is clear that reduced version can always be constructed. For
primary branchwise equivalences, they can be the same. Now, every
other branchwise equivalence, as we have seen, is created by
subsequent refinements of the primary branchwise equivalence. Thus, we
can simply skip the refinements of the branches which are not
intersected by the ray. There is one exception to this rule, if the
primary branch which contains the ray's end is boundary-refined so
that new branches accumulate at the endpoint not covered by the ray.
Then, the next adjacent branch has to be boundary-refined, too. That
is why we allowed two refinements in the definition. 

\subparagraph{The marking.}
The {\em marking} determined by the ray is an operator which changes
the branchwise equivalence on the branch which contains the end only,
so that the end gets mapped onto its conjugate point. 

This is not really a definition, since there are certainly many way in
which one could mark in this sense. We postpone the precise definition
until next section. Right now we only need to believe that the marking
operation has been defined. 

\subparagraph{Marked versions of branchwise equivalences.}
The {\em marked version} of a branchwise equivalence with respect to a
marking condition is, first, reduced accordingly to this condition.
Secondly, is marked in the sense of the previous paragraph.

\subparagraph{A brief explanation.}
Marked versions of branchwise equivalences will be used to be
pulled-back by critical branches. A careful reader has likely guessed
the meaning of the ray's end, which is at the critical value. Thus,
marking ensures that the pull-back is well-defined.

The interpretation of the direction of the ray is that it covers the
image of the folding branch. Thus, it makes no difference how the
branchwise equivalence is defined beyond the ray. However, when we
complexify the procedure, that region will have a preimage beyond the
real line. So, the idea is to make the branchwise equivalence as simple as
possible in that region, and avoid future trouble.

\paragraph{Summary.}

\begin{prop}\label{prop:11,2}
If the construction starts with a primary branchwise equivalence as described,
continues through an arbitrary number of steps, and ends with a final 
refinement, the result is a preferred
branchwise equivalence with the following uniform geometric estimates:
\begin{itemize}
\item
The ratio of any two adjacent domains is bounded.
\item
The lengths of the external branches are bounded away from $0$.
\end{itemize}
\end{prop}
\begin{proof}
We only give an outline. The argument for the first part was discussed.
To see that the second part is true, we notice that if an external branch
of the primary map is refined, the resulting external branch will never be
refined. Indeed, the only possibility that could happen is a chain boundary-
refinement. However, by our construction, the image of a folding branch is
never contained in an external branch of the primary map. Also, the 
refinement of the primary external branch must have created folding branches.
But, if there are folding branches between the tip of the folding branch
and the external branch, the chain boundary refinement will stumble on them
and will never reach the new external branch.
\end{proof}

If we can prove that this map is uniformly qs, our main theorem will follow.
To prove that the conjugacy is qs on the whole interval we can use arguments
similar to those in \cite{yours}. 

Thus, our main theorem reduces to the following:

\begin{em}
Prove that all branchwise equivalence obtained in the procedure described 
above are uniformly quasisymmetric. 
\end{em}

This is what the balance of the paper is about.

\section{Complexification of induced maps}
\subsection{Introductory remarks}
\paragraph{Why complexification?}

There are two basic reasons why we want to work with a complexified version
of our problem. First is that the critical pull-back is hard to handle
using the real variable methods only. This is an operation which involves
two maps with unbounded real distortion. True, their effects should cancel,
but there is no good way to account for that if we confine ourselves to the
real line.. On the other hand, since the
quadratic polynomial is analytic it has a null impact on quasiconformal 
distortions.

Another advantage of using quasiconformal maps was mentioned at the
end of the previous paragraph. The point is that is easier to paste 
quasiconformal maps. We simply need to check that the result is continuous
and the quasiconformal distortion is bounded. Also, quasiconformal distortion
is something which can be localized. To make this point more clear let us
consider two real quasisymmetric maps: one is the identity to the left of
zero, the other to the right. It is intuitively obvious that since their 
``quasisymmetric distortions'' are supported in different regions, the 
quasisymmetric norm of the composition should be more like a maximum than
the sum of the norms. But there is no correct and convenient way to express 
this kind of intuition other than in terms of quasiconformal 
maps. The support of quasiconformal distortion is a well-defined notion and 
the fact the distortions are supported in different regions can be used
correctly.   

\paragraph{The strategy.}
In this section, we will extend branches of induced maps to quasiconformal
mappings of the whole plane. We will require these extensions to have
special properties and, in fact, the main technical burden of our work is
going to be in that part.

In the next section, we will be able extend the branchwise equivalences first 
for all
branches more or less independently, and only on a small set around each
branch. We will show that the ``piecewise extensions'' obtained in this
way are uniformly quasiconformal where defined. 

Finally, we will show that the piecewise extensions can glued and that only
involves bounded quasiconformal distortion.

\subsection{Extensions of individual branches}
The objective of this passage is to show how a branch can be extended to the
whole plane. It is important that at this moment we regard the branch as a
separate entity. That means, we will not care about whether our extension
is consistent with other branches. 

\paragraph{Simple extensions.}
\subparagraph{Tangent extension.}
Suppose we have a monotone extendable branch $\xi$ defined on an interval $J$.
Let us extend $\xi$ on the whole line using affine maps so that the resulting
map is differentiable. By the definition, the {\bf tangent extension} of
$\xi$ is the tangent map of the function defined in this way, where tangent 
spaces are identified with vertical lines.

It should be mentioned here that the idea of tangent extensions and their 
properties were discussed in D. Sullivan's lectures in New York in the fall of
1989. 

The qc distortion of tangent extensions can be computed in an elementary way,
and the result can be expressed as follows:
\begin{fact}\label{fa:n6,1}
Consider a monotone extendable branch $\xi$ and its tangent extension 
$\ex(\xi)$. Rescale by an affine mapping so that the length of the domain of 
$\xi$ becomes $1$. Then, the conformal distortion at $(x,y)$ equals
\[ \frac{\partial_{\overline{z}}\ex(\xi)}{\partial_{z}\ex(\xi)} =
\frac{iy\cdot f''/f'(x)}{2-iy\cdot f''/f'(x)} \; .\]
\end{fact}

\paragraph{Local extensions of branches.}

\begin{defi}\label{defi:21,1}
Consider a monotone branch $\xi$. It is defined as an iterate of $f$ on its 
domain, and
the iterate of $f$ can be represented as an alternating composition of 
maps $h$ and the quadratics. The {\em local extension} of $\xi$ is
defined as the corresponding composition of maps of the plane, where 
transformations $h$ restricted the images of the domain have been replaced by 
their tangent extensions, and the quadratics have been extended analytically. 
\end{defi}

\begin{defi}\label{defi:n6,5}
Definition ~\ref{defi:21,1} is extended on folding branches as
follows. By definition ~\ref{defi:9,4}, a folding branch is a
composition of one monotone branch, a quadratic polynomial, and
another monotone branch. To obtain its local extension, we extend the
monotone branches locally, and the quadratic analytically.
\end{defi}

\subparagraph{Problems with local extensions.}
Eventually, we will cut off small pieces of local extensions around the domains
on the real line and will glue them into a global quasiconformal map. But in 
order to do that, we at least need to know that ``small pieces around the real
domains'' map to ``small pieces around the real images'' and not in a totally
weird way. Our next, highly technical section, is devoted proving the suitable 
estimates.

\subsection{Local extensions in the proximity of the real line}
\paragraph{Normalized monotone branches.}
\begin{defi}\label{defi:21,2}
Suppose that we have diffeomorphism written as a composition
\[ \xi=h_{n} \circ Q_{n} \circ\cdots\circ h_{1}\circ Q_{n} \]
where $h_{i}$ are negative Schwarzian maps, while $Q_{i}$ are quadratic 
polynomials with critical points $c_{i}$ respectively. Furthermore, we
assume that all maps are automorphisms of the
unit interval, and that the composition is defined, and still a
diffeomorphism,  on a larger
interval $J$ so that $\xi(J)=(-\epsilon,1+\epsilon)$ where $\epsilon$ is a 
constant between $0$ and $1$ to be specified soon. 

We will call this composition a {\em  normalized monotone branch}. 

Correspondingly, we can consider its local extension, in which 
maps $h_{i}$ are extended tangentially, whereas the polynomials are extended 
analytically.
\end{defi}  

We will use the notation $J_{0}:=J$ and $J_{i}=h_{i}(Q_{i}(J_{i-1}))$.
\subparagraph{Correspondence between the local extensions of branches
and local extensions of normalized branches.}

Let us consider a monotone extendable branch with a stopping rule $s$. 
Look at all $s$ images of
the domain. They can all be affinely rescaled to become the unit interval.
Thus, we get a corresponding normalized branch. Its local extension
again is related to the local extension of the branch by affine maps.   
Also, the extendability condition is satisfied with a uniform $\epsilon$.  

There is an important estimate for normalized  monotone branches which
come from monotone extendable branches in the way just described.
\begin{fact}\label{fa:n6,2}
If an abstract monotone branch$\theta$ comes from an extendable monotone
branch, then $S\theta > -K_{1}$.
\end{fact}
\begin{proof}
It is a standard fact that this follows from extendability. See
\cite{harm}, Proposition 2, for the proof of a very similar statement. 
\end{proof}

\subparagraph{The relation of standard and analytic extensions of the 
quadratics.}

If the polynomials $Q_{i}$ were extended tangentially, then 
local extensions would be 
easy to understand. So, it is natural to study a relation between 
the analytic extension of a quadratic map and the tangent extension of the
same map. We have a Lemma:
\begin{lem}\label{lem:20,1}
Let us consider a quadratic map $F$ from the unit interval into itself, with 
the critical point in $c\not\in I$. For a point $z=(x,y)$, $0<x<1$,
the $y$ coordinate of the image is the same for both the tangent and the 
analytic extension. The $x$ coordinates differ not more than 
\[ K_{1} \frac{y^{2}}{(\mbox{dist}(c,I))^{2}}\; .\]
\end{lem}
\begin{proof}
This elementary geometry.
\end{proof}

\paragraph{Perturbations of normalized monotone branches.}
If we pick a point $z=(x,y)$ and look at its image by the local extension of 
$h_{i}\circ Q_{i}$,
we discover that the $x$-coordinate of the image is not equal to
$h_{i}\circ Q_{i}(x)$. However, the discrepancy is rather small as
indicated by Lemma~\ref{lem:20,1}. Thus, to follow the $x$-coordinates
of the images of a point by consecutive maps from the composition, we
need to perturb the normalized monotone branch. This gives rise to the
following object.

\begin{defi}\label{defi:q13a,1}
If $\xi=h_{n}\circ Q_{n}\circ\cdots\circ h_{1}\circ Q_{1}$ is a
normalized local branch, its {\em perturbation} is any composition 
\[\overline{\xi}=h_{n}\circ g_{n}\circ Q_{n}\circ\cdots\circ h_{1}\circ
g_{1}\circ Q_{1}\]
where each $g_{i}$ is an orientation-preserving homography that fixes 
$Q_{i}(J_{i-1})$. 
\end{defi}

If $g_{i}$ fixes $(a,b)$, it is uniquely characterized by the number
\[ \Delta_{i}=\log\frac{(x-a)(g_{i}(x)-b)}{(g_{i}(x)-a)(b-x)}\]
which is independent of the choice $x$. 

Before we consider more closely the correspondence between local
extensions and perturbations, here is a the property of perturbations
which will be of interest to us.

\begin{lem}\label{lem:q13a,1}
Consider a perturbation $\overline{\xi}$ and a point $x$ so that 
$\overline{\xi}(x)\in [0,1]$. Then, $\overline{\xi}'(x)$ is bounded
away from $0$ by a constant. Moreover, provided
$\sum_{i=1}^{n}|\Delta_{i}|$ is sufficiently small, it is also bounded
in uniform way.
\end{lem}
\begin{proof}
Denote $J_{0}=(a,b)$.

The first claim follows immediately because the infinitesimal cross-ratio
\[ \frac{(b-a)dx}{(x-a)(b-x-dx)} \]
is increased by $\overline{\xi}$ (to see why, consult \cite{preston}).

The second claim will follow if we can show that $\overline{\xi}^{-1}(0,1)$
has uniformly large length.

Consider the sequence \[u_{i}=Q_{n-i+1}^{-1}\circ g_{n-i+1}^{-1}\circ
h_{n-i+1}^{-1} \circ \cdots\circ Q_{n}^{-1}\circ g_{n}^{-1}\circ 
h_{n}^{-1}(0)\; .\]

If $J_{n-i}=(a',b')$. We claim that
\[ |\log\frac{(0-a')(b'-u_{i})}{(u_{i}-a')(b-0)}|\leq\sum_{n-i+1}^{n}
|\Delta_{i}|\; .\]
Indeed, this is immediately seen true by induction if one keeps in
mind that $Q_{i}^{-1}\circ h_{i-1}^{-1}$
contracts the ``Poincar\'{e} metric'' on the images of $J$. The
Poincar\'{e} metric on an interval is the conformally invariant
metric on the disc whose diameter is the interval. It is classically
known that it can be represented as the logarithm of some cross-ratio,
thus it is expanded by negative Schwarzian maps and contracted by
positive Schwarzian maps.
      
The same argument can be applied to the preimages of $1$. Since the 
Poincar\'{e} length of the interval $(0,1)$ inside $J_{n}$ is 
more than $-2\log\epsilon$, it is enough for
$\sum_{i=1}^{n}|\Delta_{i}|$ to be less than $-\log\epsilon/2$ in
order to ensure that the Poincar\'{e} length of
$\overline{\xi}^{-1}(0,1)$ is definite.
\end{proof}

\subparagraph{Orbits by the local extension and perturbations.}
We consider a sequence $z_{0}=(x_{0},y_{0})$ and 
\[z_{i}=(x_{i},y_{i}):=h_{i}\circ Q_{i}(z_{i-1})\; .\]  
If also $0<x_{i}<1$ for every $i$ between $0$ and $n-1$ inclusively,
we call this sequence an {\em orbit}.

With every orbit we can associate the unique perturbation which satisfies
\[ x_{i}=h_{i}\circ g_{i}\circ Q_{i}(x_{i-1}) \] for all $i$ from the
relevant range. One can check that $\Delta_{i}$ is bounded
proportionally to the discrepancy between the tangent and analytic
extension of $Q_{i}$ at $z_{i-1}$, thus by $K_{1}\frac{y^{2}_{i-1}}
{\dist([0,1],c_{i})}$ according to Lemma~\ref{lem:20,1}. Therefore,
the crucial sum $\sum_{i=1}^{n}|\Delta_{i}|$ can be estimated by
\[\sum_{i=1}^{n}|\Delta_{i}|\leq
K_{1}(\max\{y_{i:0\leq i<n}\})^{2}\sum_{i=0}^{n-1}
\frac{1}{(\dist([0,1],c_{i}))^{2}}\; .\]

On the other hand, we can calculate that for any $x\in (0,1)$
\[ S\xi(x) \leq
-\sum_{i=0}^{n-1}\frac{K_{2}}{(\dist([0,1],c_{i}))^{2}}\; .\]

Finally, as a consequence of Fact~\ref{fa:n6,2} we obtain
\[\sum_{i=1}^{n}|\Delta_{i}|\leq K_{3}(\max\{y_{i}:0\leq i<n\})^{2}\; .\]  

That the perturbation contains some information about the orbit is
evident from our next lemma.
 
\begin{lem}\label{lem:3,14}
Consider an orbit $z_{i}$ and the corresponding perturbation
$\overline{\xi}$. We have an estimate
\[|\log\frac{y_{n}}{y_{0}}-\log\overline{\xi}'(x_{0})  
\leq -\sum_{i=0}^{n-1}\frac{K_{2}}{(\dist([0,1],c_{i}))^{2}}\; .\]
\end{lem}
\begin{proof}
Observe that 
\[ y_{i}/y_{i-1}= \frac{(h_{i}\circ g_{i}\circ
Q_{i})'(x_{i-1})}{g'_{i}(Q_{i}(x_{i-1}))} \; .\]
So, the difference is bounded by the sum of logarithms of $g'_{i}$ at
the appropriate points.
However, $|\log g'_{i}(Q_{i}(x_{i-1}))|$ can be estimated according to
Lemma~\ref{lem:20,1}, and then the sum can be bounded using
Fact~\ref{fa:n6,2}. 
\end{proof}

For an orbit, let $Y$ denote the maximum of $y_{i}$ with $0\leq i<n$.
The essence of our results obtained so far is in the following lemma:
\begin{lem}\label{lem:q13p,1}
Provided $Y<K_{1}$, 
\[|\log y_{n}/y_{0}|\leq K_{2}\; .\]
\end{lem}
\begin{proof}
This is simply a summary of Lemmas~\ref{lem:q13a,1} and
~\ref{lem:3,14} together with our estimate of $\sum_{i=0}^{n-1}\Delta_{i}$. 
\end{proof}

\paragraph{The basic result.}
We now have all necessary tools to quickly prove the main result of
this section.

For a point $z=(x,y)$ with $0<x<1$, let the {\bf height} of $z$ mean
$y/\min(x,1-x)$.

\begin{prop}\label{prop:11,1}
If the height of some $z$ is less than $K_{1}$, a uniform constant, and 
$z_{n}=(x_{n},y_{n})$ is the image of $z$ under the local 
extension of a normalized monotone branch, then $0<x'<1$ and the
ratio of heights of $z_{n}$ and $z$, as well
as $y_{n}/y$ is uniformly bounded from both sides.
\end{prop}
\begin{proof}
The condition for the ratio $y_{n}/y$ is similar to the claim of
Lemma~\ref{lem:q13p,1}, but there are two things that we need to
check. First of all, we need to show that $z$ defines an orbit, that
is $x_{i}$ is always between $0$ and $1$. Secondly, we need to bound
$Y$ in terms of $y_{0}$.  

\subparagraph{Working to eliminate $Y$.}
In the following two lemmas we simply assume that $0<x_{i}<1$ for
$0\leq i<n$.

\begin{lem}\label{lem:20,4}
For a suitable constant $K_{1}$, if $y_{0}<K_{1}$, then $y_{n}< K_{2}y_{0}$.
\end{lem}
\begin{proof}
This is a simple corollary to Lemma ~\ref{lem:q13p,1}. The trick is to
look at  the first $i$ for which $y_{i}>\exp(K_{2,L.~\ref{lem:q13p,1}})y_{0}$. 
Provided $y_{0}$ is small, Lemma~\ref{lem:q13p,1} applied to the
composition cut off at $n:=i$ gives a contradiction.
\end{proof} 

\begin{lem}\label{lem:20,5}
Provided $y_{0}<K_{1}$, 
\[|\log y_{n}/y_{0}|\leq K_{2}\; .\]
\end{lem}
\begin{proof}
Follows immediately from Lemmas~\ref{lem:q13p,1} and ~\ref{lem:20,4}.
\end{proof}

\subparagraph{Working to prove that $z_{i}$ is an orbit.}

Next, we have to investigate how $x_{i}$ depends on $x_{0}$.
Fortunately, this is reduced to a one-dimensional question about perturbations.
We introduce a family of perturbations $\overline{\xi}_{t}$. If
$\Delta_{i}$ the difference between the $x$-coordinate of $Q_{i}(z_{i-1})$
and $Q_{i}(x_{i})$, then $\overline{\xi}_{t}$ corresponds to the
sequence $t\Delta_{i}$. Hence, $\overline{\xi}_{0}=\xi$ while 
$\overline{\xi}_{1}=\overline{\xi}$. We also get sequences $x_{i}(t)$
whose definition is natural. We further denote
\[ \overline{\xi}_{t}^{i} := h_{n}\circ g_{n,t}\circ
Q_{n}\circ\cdots\circ h_{i}\circ g_{i,t}\circ Q_{i}\]   
 
\begin{fact}\label{fa:3,14}
\[ \frac{dx_{n}(t)}{dt} = \sum_{i=1}^{n}
\frac{d\overline{\xi}_{t}^{i}}{dx}(x_{i}(t)\cdot
h'_{i-1}(h^{-1}_{i-1}(x_{i}(t))\Delta_{i} \; .\] 
\end{fact}
\begin{proof}
Elementary.
\end{proof}

Our next lemma is analogous to Lemma~\ref{lem:q13p,1} except that
the $x$-coordinate is now involved.
\begin{lem}\label{lem:20,6}
For any $K_{1}>0$, a $K_{2}>0$ can be chosen so that
if the height of $z$ less than $K_{2}$ and $0<x_{i}(t)<1$ for all  
$0\leq i<n$ and $0\leq t\leq \tau\leq 1$, then 
$|x_{n}(\tau)-x_{n}(0)| < K_{1}\min(x_{0},1-x_{0})$.
\end{lem}
\begin{proof}
If the height of $z$ is small, Lemma~\ref{lem:20,5} applies. This
allows us to estimate $\sum_{i=0}^{n}\Delta_{i}$ by
$K_{2}K_{3}\min(x_{0},1-x_{0})^{2}$. Subsequently,,
Lemma~\ref{lem:q13a,1} can be applied to estimate the derivatives in
the formula of
Fact~\ref{fa:3,14} by constants. Then, that formula immediately yields
the claim. 
\end{proof}

Finally, we notice that  Lemma~\ref{lem:20,6} remains true even when
the assumption of $x_{i}(t)$ being between $0$ and $1$ has been
removed. To this end, we note that
$\min(x_{n}(0),1-x_{n}(0))/\min(x_{0}(0),1-x_{0}(0))$ is bounded away
from $0$. That follows directly from extendability of the monotone
branch. So, $K_{1,L.~\ref{lem:20,6}}$ can be chosen so as to ensure
that $|x_{n}(\tau)-x_{n}(0)| < \min(x_{n}(0),1-x_{n}(0))$. Then, we
consider an arbitrary $z$ which satisfies other assumptions of 
Lemma~\ref{lem:20,6} and look for the lowest $\tau\leq 1$ and lowest
$i$ so that $x_{i}(\tau)\notin (0,1)$. Then Lemma~\ref{lem:20,6}
applied to the composition cut off at $n:=i$ gives a contradiction.

Proposition~\ref{prop:11,1} follows directly.
\end{proof}

\subsection{Global extensions}
\paragraph{Diamond neighborhoods.}
\begin{defi}\label{defi:n6,3}
For an interval, its {\em diamond neighborhood} of size 
$a$ is defined to be an open quadrilateral bounded by the set of points with 
height $a$.
\end{defi}

We have a lemma which is a simple corollary to Proposition
~\ref{prop:11,2}:
\begin{lem}\label{lem:n6,3}
Fix an $\alpha\leq K_{1,P.\ref{prop:11,1}}$. Then,  there is a
positive function $\Gamma$ so that, for any extendable branch 
$\xi=h_{2}\circ Q\circ h_{1}$, the
local extension of $\xi$
maps the size $\Gamma(\alpha)$ diamond neighborhood of the domain of $\xi$   
into the size $\alpha$ diamond neighborhood of the image of $h_{2}$. 
\end{lem}
\begin{proof}
First, assume that $\xi=h_{1}$, i.e. $\xi$ is monotone. Then, the
claim follows immediately from Proposition ~\ref{prop:11,1}. Next,
consider the case of $\xi=Q\circ h_{1}$. We see in an elementary way
that  $Q$ will not extend a
diamond neighborhood of bounded size to much.
To conclude the argument, we apply Proposition~\ref{prop:11,1} to $h_{2}$. 
\end{proof}

\begin{defi}\label{defi:n6,4}
The diamond neighborhood of size
$\min(1/2,\Gamma(K_{1,P.\ref{prop:11,1}}))$, in the
notations of Lemma ~\ref{lem:n6,3} of the domain of any branch will be
called the {\em large diamond neighborhood} of that branch.
\end{defi}

\begin{defi}\label{defi:q21a,2}
A diamond-like neighborhood of an interval on the real line is any
open neighborhood of the interior of the interval. The size of a
diamond-like neighborhood is defined to the size of the largest
diamond neighborhood of that interval contained in the diamond-like
neighborhood. 
\end{defi}

\paragraph{Introduction of global extensions.}
\subparagraph{Postulates.}
We will construct a new type of complex extension of the branch,
called a {\em global extension}. The idea is to make the same as the
local extension close to the real domain of the branch, but change it
far  from that domain in order to make glueing possible with
extensions of other branches. 

Thus, formally, the global extension of a branch $\xi$ will satisfy
these postulates:
\begin{enumerate}
\item
The global extension is the same as the local extension on the large
diamond neighborhood of $\xi$.    
\item
If the domain of $\xi$ is $(a,b)$, the global extension is affine
outside of the rectangle with vertices at $a+(b-a)i$, $b+(b-a)i$,
$b-(b-a)i$, and $a-(b-a)i$.
\item
The global extension is uniformly qc function.
\end{enumerate}

\subparagraph{The construction.}
We will sketch the construction of global extensions from local
extensions. First, we consider the case of a monotone branch and its
corresponding normalized monotone branch. We choose approximate maps so
that the extension remains local on the diamonds of size
$K_{1,P.\ref{prop:11,1}}$ around each intermediate image of the domain,
while it is tangent outside of the diamond neighborhood of unit height.
We leave without a proof that such adjusting map can be constructed
with complex distortion of the order of $\dist((0,1), c_{i})^{-2}$.  
This, in view of Fact ~\ref{fa:n6,2}, ensures that the complex distortion
of the composition will be uniformly bounded. 

Having thus constructed a map which is a tangent extension outside a
diamond of unit height, it is easy to build a
global extension, and we leave it to the reader. 

In the case of $\xi=h_{2}\circ Q\circ h_{1}$, we can build global
extensions of both monotone branches $h_{1}$ and $h_{2}$, and the
suitable extension of $Q$ can be constructed in an elementary way.

\subparagraph{A convention.}
In the future will work mostly with global extensions. So, if we say
``an extension of the branch'' we mean the global extension. In
formulas, the global extension of $\xi$ will appear as $\ex(\xi)$. 
\paragraph{Better qc estimates for global extensions.}
We will show a lemma about the qc distortion inside large diamond
neighborhoods. 

\begin{lem}\label{lem:n6,4}
Consider a monotone extendable branch. Let $D$ mean the total length
of all intermediate images of the domain of this branch. The complex
distortion of the extension of this branch at a point $z$ inside the large
neighborhood is bounded by a constant multiplied by $D$ and by the
ratio of the distance from $z$ to the line to the length of the domain
of the branch.
\end{lem} 
\begin{proof}
Consider the corresponding normalized branch. The key estimate is Fact 
~\ref{fa:n6,1}. The complex distortion of the composition is bounded
by the sum of complex distortions of maps $h_{i}$ along the orbit of $z$. 

According to Fact ~\ref{fa:n6,1}, each contribution is bounded
proportionally to $y_{i}h''_{i}/h'_{i}(x_{i})$. 
The quantity $y_{i}$ is roughly constant according to
Proposition ~\ref{prop:11,1}, and so it is comparable to the ratio of
distance from $z$ to the line by the length of the domain.

To estimate estimate the ``nonlinearity ratio'' $h''_{i}/h'_{i}$ we
need to remember that $h_{i}$ is
just an affine rescaling of $h$ restricted to the $i$-th image of the
domain. Clearly, $h''/h'$ is uniformly bounded, and affine rescaling
multiplies the nonlinearity ratio by the derivative of the rescaling
map, which is equal to the length of the domain of $h_{i}$ in our
case.
\end{proof}

Lemma ~\ref{lem:n6,4} is useful provided that we can give a good
estimate of $D$. Fortunately, it is so. Since there is a lot of
expansion in our inducing construction, we will be able to show that
in interesting cases the $D$ is just proportional to the length of the
image of the branch. 

Lemma ~\ref{lem:n6,4} also has obvious consequences for folding
branches, since they can be written as
\[ h_{1}\circ Q\circ h_{2}\] 
with $h_{1}, h_{2}$ monotone and $Q$ analytic. Lemma~\ref{lem:n6,4}
may be applied to both monotone branches separately.

\section{Complexified branchwise equivalences}\label{sec:3}
\subsection{Basic properties and constructions}
\paragraph{Primary branchwise equivalences.}
In the section on real branchwise equivalences we described how to
construct primary branchwise equivalences on the real line. The main
objective of this section is to describe the process of their
complexification. However, even before we do that, we wish to discuss
an important technical principle of the construction.
\subparagraph{Arbitrary fineness principle.}
Roughly speaking, we may assume that the longest domain of the primary
induced map is as short as it suits us. More precisely, we claim 
\begin{fact}\label{fa:n6,4} 
For any $\alpha>0$ primary induced maps can be constructed so that the
lengths of the branches do not exceed $\alpha$, and the qs norm of the
corresponding branchwise equivalence is bounded by $\Gamma(\alpha)$ where
$K$ is the function of $\alpha$ only.
\end{fact}
\begin{proof}
We have shown how to construct uniformly qs primary branchwise
equivalences without the fineness requirement. If we want a finer
branchwise equivalence, we need to refine by pull-backs. If we
simultaneously refine all the branches by pull-backs, the lengths of
the domains will decrease by a fixed factor. So, only a finite number
of such steps will be needed to attain the specified fineness. Each
step, however, increases the qs norm by a bounded amount, which can be
seen by arguments analogous to those given in the Addendum. Also, a
less tedious complex way will shown later to see that.   
\end{proof}

In the future, we will often assume that ``the primary branchwise
equivalence is sufficiently fine.'' That, in effect, means that the
$\alpha$ which occurs in Fact ~\ref{fa:n6,4} will be chosen many
times. However, since this paper will hopefully end up having a finite
length, a positive minimum will still exist. 

\subparagraph{Primary extensions of branchwise equivalences.}
Here, we list list the desired properties of complex primary branchwise
equivalences. We assume that the domains of both induced maps have
been normalized by affine maps to become $[-1,1]$.

\begin{enumerate}
\item
On the real line, the map is a branchwise equivalence as
described in the Branchwise Equivalences section. 
\item
The map is is the identity outside of the ball of radius $3$ centered
at $0$ on the plane. 
\item 
on the diamond neighborhood of size $1$ around each
branch, the map is affine. 
\item 
The map is quasiconformal, and its qc norm can be bounded by a uniform
function of the qs norm of the underlying real branchwise equivalence.
\end{enumerate}

The reader may note that these postulates leave us with the freedom to
define the map inside the diamonds of size $1$, as long as the map is
qc and and affine on the boundary. This freedom will be used for marking.

\subparagraph{The complexification lemma.}
\begin{lem}\label{lem:n6,5}
Suppose that a qs branchwise equivalence $\upsilon$ is given between two
induced maps or their boundary-refined derivatives which satisfy the
property that the ratio of any two adjacent branches is uniformly bounded.  
If $\upsilon$ is the identity outside of $(-3,3)$, then its primary
extension can be constructed.
\end{lem}
\begin{proof}
Draw a circle of radius $3$ from $0$. Consider a Jordan curve $\Gamma$
which consists of the upper half of the circle, the upper halves of
the boundaries of size $1$ diamond neighborhoods of branches, and the
connecting pieces of the real line from $-3$ to $-1$ and from $1$ to
$3$.

An analogous curve $\hat{\Gamma}$ exists in the phase space of $\hf$. 
Define a homeomorphism $G$ of $\Gamma$ onto $\hat{\Gamma}$ which is the
identity on the circle and the pieces of the
real line, and affine on the boundary of each diamond. The lemma will
clearly follow if we prove that this homeomorphism can be extended to
the region encompassed by $\Gamma$.

We notice that $\Gamma$ is a uniform quasicircle. An easy way to
see that is by noticing that the three point property of Ahlfors
(see \cite{ahl}).  The homeomorphism $G$ is also quasisymmetric in
the sense that if distances $|x-y|$ and $|x-z|$ are comparable, so are
the corresponding distances in the images. Moreover, the
quasisymmetric norm of $G$ is bounded in terms of the qs norm of the
map on the real line.

Next, $\Gamma$ and $\hat{\Gamma}$ can be
uniformized to the round circle by a qc map, and the counterpart of
$G$ on the boundary quasisymmetric. Since it can be extended in the
classical way (see \cite{abu}), and pulled back to the inside of
$\Gamma$, the lemma follows.
\end{proof}

Thus, primary extensions of primary branchwise equivalences exist, and
are uniformly quasiconformal. In the future, when we talk of complex
primary branchwise equivalences, we mean exactly primary extensions of
primary branchwise equivalences.

\paragraph{Admissible extensions of branchwise equivalences.}
Now we will define the class of complex extensions with desirable
properties, which will be called {\em admissible extensions.}

\subparagraph{Postulates.}
A complex extension $\Upsilon$ of a branchwise equivalence is
admissible if it satisfies the conditions listed below.
\begin{enumerate}
\item
Admissible extensions are the identity outside of the ball $B(0,3)$.
\item
They are quasiconformal.
\item
For every branch $\xi$ with stopping time $s$, we can choose an integer
$\overline{s}\leq s$ so that the corresponding settled branch
$\overline{\xi}$ allows us to represent $\Upsilon$ by the formula:
\[ \Upsilon=\hf^{-\overline{s}}\circ A\circ f^{\overline{s}}.\]  
The formula is supposed to valid on some maximal diamond-like neighborhood of
the domain of $\xi$, which we will call a {\em close neighborhood} of
that domain and denote with $D(\xi)$. The map $A$ is supposed to be
quasiconformal and affine outside of the image of $D(\xi)$ by
$\overline{\xi}$. Moreover, if $\xi$ is monotone, $A$ is simply affine.
\end{enumerate}

The primary branchwise equivalences are admissible with
$\overline{s}=0$ on all branches.

\subparagraph{The norm of admissible extensions.}
By the norm of an admissible branchwise equivalence we will mean the
maximum of its quasiconformal norm and the reciprocals of the sizes of
its close neighborhoods.
 
\subsection{Complex pull-backs.}
\paragraph{Complex marking.}
We will show how to mark an admissible complex branchwise equivalence. 
The idea is to use the last property of admissible equivalences and
change $A$ inside the image of $D(\xi)$ only. Remember, that we only
mark monotone branches, thus $A$ is affine. By Proposition~\ref{prop:11,1}, 
the image of $D(\xi)$ by $\overline{\xi}$ contains a diamond $D'$ on
size comparable to the size of $D(\xi)$. So, we change $A$ only inside
$D'$ to make it linear-fractional inside the diamond of half the size
of $D'$.   

It is easy to observe that the quasiconformal norm of such a map
depends on two estimates: the nonlinearity of the linear-fractional
map on the real line, and the smallness of the size of $D'$. 

The first bound always holds:

\begin{fact}\label{fa:n7,1}
Consider any point contained in the domain of a monotone branch, and
its image by conjugacy. Uniformize this domain, and the corresponding
domain  of $\hf$ by the unit interval using affine maps. Then, the
Poincar\'{e} distance between the point and its image is uniformly bounded. 
\end{fact}
\begin{proof}
It is enough to consider the conjugacy near the boundary of a monotone
branch.

As monotone branches have uniformly bounded distortion,
the question reduces the primary boundary-refined branchwise
equivalence. 
By construction, the dynamically defined objects scale exponentially
near $q$ and $-q$ and the exponential rate depends on the eigenvalue
at $q$.  This implies that the conjugacy moves points around only by
bounded Poincar\'{e} distances. See \cite{mike} for a detailed
analysis in a closely related case.
\end{proof}

The second estimate depends directly on the norm of the branchwise
equivalence being refined. We will address that issue later.

\paragraph{Simple pull-backs.}
We assume that an admissible branchwise equivalence $\Upsilon_{1}$ is
given. Also, a branch $\xi$ is chosen in the corresponding induced
map. Also, another admissible branchwise equivalence $\Upsilon_{2}$ is
given 
so that if $\xi$ is folding, its critical value is a monotone domain of
$\Upsilon_{2}$ and $\Upsilon_{2}$ is marked by the ray which starts at
the critical value and covers the image of the interval by $\xi$.  

We want to refine $\xi$ by pulling-back $\Upsilon_{2}$. This will be a
multi-step process.

We now show the first step in which we construct the correct map on
and close to the domain of $\xi$, but do not concern ourselves with
how this map matches the global $\Upsilon_{1}$.

\subparagraph{Extension of $\Upsilon_{1}$ on the close neighborhood.}
Since $\Upsilon_{1}$ is admissible, on $D(\xi)$ it can be represented
as
\[ \hat{f}^{-\overline{s}}\circ A\circ f^{\overline{s}}\; . \]

On the real line, the pull-back of $\Upsilon_{2}$, assuming
appropriate marking, can be written as
\[ \hat{\xi}^{-1} \circ\Upsilon_{2} \circ f^{s-\overline{s}}\circ
A^{-1}\circ\hf^{\overline{s}}\circ\Upsilon_{1} = G\circ\Upsilon_{1} \]
where $G$ is just a new notation for the complicated composition in
front of $\Upsilon_{1}$. Moreover, this formula makes sense on the
whole plane, and on $D(\xi)$ it gives the same as
\[ \hat{\xi}^{-1}\circ\Upsilon_{2}\circ\xi \;.\]

The map $G$ will be called the {\em simple pull-back} of
$\Upsilon_{2}$ by $\xi$. 

The map $G$ is the identity beyond the ball
centered at the midpoint of the domain of $\xi$, of radius three times
the length of the domain. To see that, we notice that the composition 
\[f^{s-\overline{s}}\circ A^{-1}\circ\hf^{\overline{s}} \]
is affine beyond the preimage of that ball by definition of global
extensions and the requirement imposed on $A$ by admissibility of
$\Upsilon_{1}$. Similarly,
\[\hat{\xi}^{-1}\]
is affine beyond that ball. Finally, $\Upsilon_{2}$ is the identity
except on the ball, again by admissibility.
 
\subparagraph{Where we stand with the construction.}
To define the new branchwise equivalence by $G\circ\Upsilon_{1}$ is
not quite a good idea since $G$ is not required to be the identity on
the line everywhere beyond the domain of $\xi$, and, in fact, could
not be for a boundary-refined $\Upsilon_{2}$. So, we will have to
change the simple pull-back a little bit to take care of this problem.

In our construction the final result of the refinement will still be
$G\circ\Upsilon_{1}$ on the large diamond neighborhood of $\xi$. 
So, the resulting map is going to be admissible.

Finally, it should be pointed out that $G$ ``lives'' in the extension
of the phase space of $\hf$ unlike branchwise equivalences which go
from the phase space of $f$ to the phase space of $\hf$.
\paragraph{Hexagonal extensions.}

We will now change the simple pull-back to get another map called the 
``hexagonal extension''. The hexagonal extension is not a
homeomorphism of the entire plane, but only of some hexagon around the
domain being refined. On the other hand, it is easy to extend by the
identity if $\Upsilon_{2}$ is not boundary-refined, or glued with an
analogous map around the adjacent domain of the chain refinement.
\subparagraph{Angular squeezing.}
We start with a simple extension $G$.

We will describe the procedure in polar coordinates around $0$.
We define a map 
\[{\cal S}_{0} : (-\pi,\pi) \rightarrow (-\pi/3,\pi/3)\]
which keeps everything inside the arc $-\pi/4,\pi/4$ fixed, and squeezes
the sectors $(-\pi,-\pi/4)$ and $(\pi/4,\pi)$ diffeomorphically into
$-\pi/2+0.001,-\pi/4$ and $\pi/4,\pi/2-0.001$ respectively. If we assume the distance
from $0$ is unchanged, this defines through polar coordinates the map also
denoted by ${\cal S}_{0}$ which is quasiconformal and can be extended to a 
multivalued function through the negative numbers.

An analogous procedure can be carried out around $1$ and the resulting 
map is to be denoted ${\cal S}_{1}$.
Then, we may consider the map
\[\rho_{1}={\cal S}_{0}\circ {\cal S}_{1}\circ G\circ 
{\cal S}_{1}^{-1}\circ {\cal S}_{0}^{-1}\; .\]

\subparagraph{Vertical squeezing.}

We would like to modify the map in such a way that all the
above listed properties remain true and the last one holds with $|Im z|<1$.

This can be easily done by introducing a diffeomorphism $\cal V$
which is the identity inside the strip $-0.5 < \Im z < 0.5$.
The strip $-1 < \Im z < 1$ which gets mapped onto $-2 < \Im < 2$ and 
outside that strip the map is a shift by $1$ vertically. Obviously,
a quasiconformal map with these properties exists.

So then we may consider
\[\rho_{2} := {\cal V}^{-1}\circ\rho_{1}\circ{\cal V}\]
which does what we wanted.
\subparagraph{A definition and comments on hexagonal extensions.}
A hexagonal extension, denoted by $G_{h}$ is defined to be 
$\rho_{2}$ restricted to the hexagon with vertices
$0,0.01-i,0.99-i,1,0.99+i,0.01+i$. It is easily verified that indeed
this hexagon is inside the domain of $\rho_{2}$. 

The hexagonal extension $G_{h}$ has a number of properties which will
be important for us
and can be verified straightforwardly:
\begin{itemize}
\item
It is the same as $G$ inside the diamond with vertices
$0,0.5-0.5i,1,0.5+0.5i$, in particular on the large diamond
neighborhood of the domain of $\hat{\xi}$.
\item
The map can be continuously extended and then the boundary of its domain is
mapped onto itself.
\item
The map is quasiconformal.
\item 
It is the identity outside the region $|Im z| < 2$.
\item
$G_{h}$ is the identity on the top and bottom edges, and also on its
entire boundary if $\Upsilon_{2}$ was not boundary-refined.
\end{itemize}

\paragraph{Pasting the neighbors in the chain boundary-refinement.}
We are ready to describe how to complexify a chain boundary refinement.

First, can construct hexagonal extensions for all members of the chain.
 The next thing to do is to glue the neighbors together.
Namely, we assume that one interval of the chain is $(0,1)$ and its neighbor 
is $(-a,0)$.  All this looks as follows:

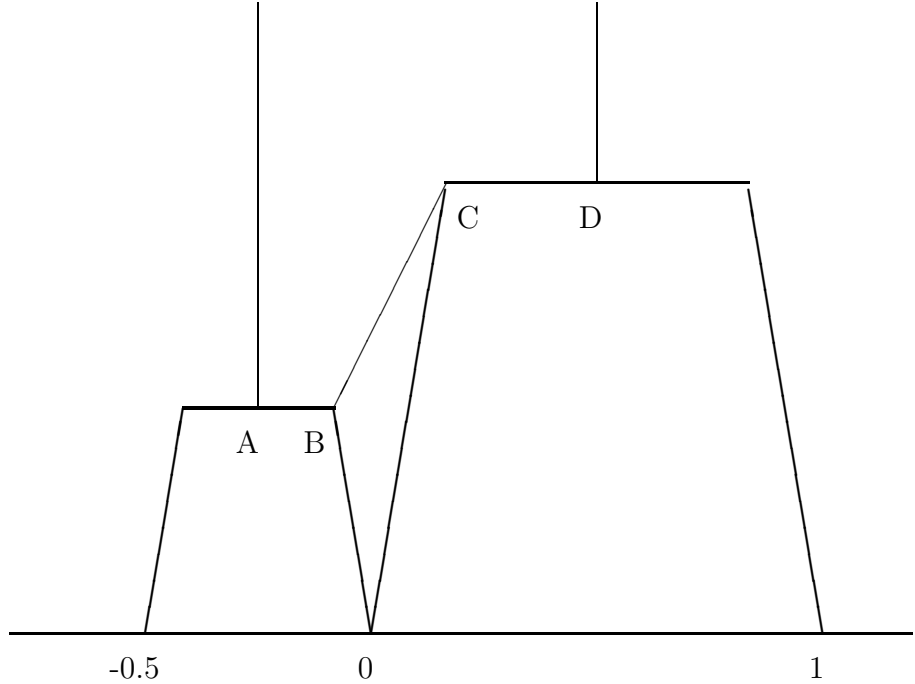
\setlength{\unitlength}{0.6mm}
\begin{picture}(200,170)(-80,-10)
\thicklines
\put(-80,10){\line(6,0){200}}
\put(0,10){\line(1,6){16.4}}
\put(0,10){\line(-1,6){8.3}}
\put(100,10){\line(-1,6){16.4}}  
\put(16.25,110){\line(6,0){67.7}}
\put(-50,10){\line(1,6){8.3}}
\put(-41.88,60){\line(6,0){33.8}}
\thinlines
\put(16.75,110){\line(-1,-2){24.93}}
\put(-25,60){\line(0,6){90}}
\put(50,110){\line(0,6){40}}
\put(-30,50){A}
\put(-15,50){B}
\put(19,100){C}
\put(46,100){D}
\put(-58,0){-0.5}
\put(-3,0){0}
\put(97,0){1}
\end{picture}

The drawing shows only half of the actual picture, which is symmetric with
respect to the real axis. The maps originally are only defined
within the two hexagons. We will define their {\bf glueing map} which  will 
extend their set-theoretical
sum. In addition to the union of the hexagons, the glueing map is defined 
in the whole infinite strip between the lines $\Re =-0.25$ and $\Im = 0.5$.
Moreover, it is assumed that the glueing map is the identity everywhere above
the line joining points $A,B,C,D$. Of course, it is also symmetric with
respect to the real axis. Finally, it is uniformly quasiconformal.\footnote{
The word ``uniformly'' meaning that the bound for the QC norm should not 
depend on a particular choice of branches, only on two maps $f$ and $\hf$.}  

\begin{prop}\label{prop8,1}
A glueing map with these desired properties can always be constructed.
\end{prop}
\begin{proof}
Since the lengths of any two adjacent branches in our
construction are comparable within uniform constants (see
Proposition~\ref{prop:11,2}), the map on $B0C$  is 
intrinsically quasisymmetric. We leave it to the reader to
complete the proof.   
\end{proof}

\subparagraph{The quasiconformal implementation of the chain boundary 
refinement.}
To implement the chain boundary refinement we
glue  together all neighbors using glueing maps and then take a set-theoretical
sum of all glueing maps. This map is then extended by
the identity on the whole complex plane. The complexified version of the 
chain boundary refinement is then composition with that map, called the 
{\em refining map}.

If the chain ends at a folding branch, we can put $G_{h}$ equal
to the identity on this folding branch, glue it with the last hexagon of the 
chain, and extend by the identity beyond the folding branch.  

\paragraph{Summary.}
We can describe the complex realization of two main pull-back
operations: the pull-back of on a single branch, and the chain
pull-back. In both cases, the new
branchwise equivalence $\Upsilon_{3}$ can be written as 
$\overline{G_{h}}\circ\Upsilon_{1}$
where $\overline{G_{h}}$ is called the refining map. The construction
of the refining map is more complicated in the case of a chain
refinement and has been described above. If only one branch is refined
by pulling back a non-boundary-refined $\Upsilon_{2}$, then the
refining map is just the extension of the hexagonal extension by the
identity outside of the hexagon. 

\subparagraph{Regions of the refining map.}
Here, we assume that inside every close neighborhood a diamond neighborhood
has been chosen, called {\em small diamond}. We will later explain how to
specify small diamonds.

Given a refining map, we can split the plane into three regions:
\begin{itemize}
\item
The pull-back region which is the union of images by $\Upsilon_{1}$ of
the close 
neighborhoods of the branches being refined.
\item
The trivial region which consists of all points whose whole neighborhoods are
fixed by the map.
\item
The glueing region which is the rest.
\end{itemize}

\subparagraph{The push-forward map.}

On the component of the pull-back region around the domain of a branch
whose stopping time is $s$, we have the push-forward
map defined simply as $\ex(f^{s})$. The image should be thought to belong   
to the domain of $\Upsilon_{2}$. Thus, formally the push-forward
``map'' is a pair: the map itself, and $\Upsilon_{2}$.

\subsection{Filling-in.}
\paragraph{The structure of the complex filling-in.}
In the Branchwise Equivalences section we described the filling-in on
the level of real branchwise equivalences as a limit of the sequence
of critical pull-backs. This allows for an immediate extension of the
procedure, since we have already defined complex realizations of
critical pull-backs. However, certain questions emerge because of the
infinite nature of this process. The most important thing that we need
to prove is that the limit exists. Also, there is an issue of whether
the limit map is going to be quasiconformal or even a homeomorphism.
We prove a lemma which immediately implies the existence of a limit.
The estimates will not be tackled until the next section.

\subparagraph{Formal complex description.}
We are given an admissible branchwise equivalence $\Upsilon$ marked by
its own critical value so that the ray covers the image of the
interval by the folding branch. We assume that, if necessary, some
branches of $\Upsilon$ have been boundary-refined and, therefore, all
monotone branches obtained as a result of the critical pull-back will
be uniformly extendable. 

We build the sequence of branchwise equivalences defined inductively
as follows:
\begin{enumerate}
\item
\[ \Upsilon_{0}:=\Upsilon\; .\]
\item
$\Upsilon_{i+1}$ is the result of a simultaneous critical pull-back of
$\Upsilon_{i}$ onto all critical branches whose domains are covered by
the marking ray.
\end{enumerate}

\subparagraph{A fineness requirement for $\Upsilon$.}
Consider two quantities: $\alpha$, the length of the longest domain of a branch
of $\Upsilon$, and $\beta$, the infimum of the sizes of small diamonds around 
the branches of $\Upsilon$. A fineness requirement  
is that the ratio $\frac{\beta}{\alpha}$ be sufficiently large. How
large we need will specified in the proof of the next lemma.

\paragraph{The existence of a limit.}
\begin{lem}\label{lem:807p,1}
Let $\Upsilon$ satisfy a suitable fineness requirement, to be
specified on the course of the proof. Then the sequence $\Upsilon_{i}$
converges everywhere. Moreover, for every point $z$ not in the real
line, its image $\Upsilon_{i}(z)$ stabilizes after a finite number of
steps.
\end{lem}       
\begin{proof}
We observe that all $\Upsilon_{i}$ are the same as
$\Upsilon_{0}$ except on some diamond neighborhood of the interval
$[-1,1]$ whose size is proportional to the length of the longest
folding branch of $\Upsilon$. Indeed, each $\Upsilon_{i+1}$ is
obtained by the critical pull-back onto the folding branches of
$\Upsilon$. However, the refining map for each pull-back is the
identity except on the hexagon of the diameter comparable to the
domain length.

From these two observations, we infer that a fineness requirement can
be chosen so that the image of the push-forward map for any branch of 
$\Upsilon$ being refined contains the region in which $\Upsilon_{i}$
and $\Upsilon_{0}$ differ. The possibility of that follows from
Proposition ~\ref{prop:11,1}.

Then, it follows immediately that $\Upsilon_{i+1}$ is the same as
$\Upsilon_{1}$ except on the preimage of the pull-back region of the
refining map. Since the map being refined is always $\Upsilon$, the
push-forward is fixed. Thus, we see that if $z$ is in the domain of
the push-forward map, the sequence $\Upsilon_{i}(z)$ stabilizes if and
only if the sequence $\Upsilon_{i}(\rho(z))$ stabilizes where $\rho$
stands for the push-forward map. Thus, the lemma will be proven if we
show that $\rho$ can only have finitely many iterates for any $z$ not
on the real line.    

This is so because, again if $\Upsilon$ is sufficiently fine, we
see that the distance of $\rho(z)$ from the real line is either larger
by a fixed constant than the distance from from $z$ to the line, or
$\rho(z)$ is no longer in the domain of $\rho$. To this end, we have
to examine the folding branch whose small diamond contains $z$. It can
be written as $h_{1}\circ Q\circ h_{2}$. The monotone branches $h_{1}$
and $h_{2}$ roughly preserve distances from the real line relative to
the domain and the image. Thus, since the domain is small and the
image is the whole interval $[-1,1]$, the distance indeed grows unless
$Q$ decreases it a lot. That can only happen if $h_{1}(z)$ is close to
the imaginary line. But then $\rho(z)$ will be close to the line, but
its projection onto the line will be beyond the marking ray. Thus,
$\rho(z)$ will be in the region where $\Upsilon_{i}$ is no different
from $\Upsilon_{0}$.
\end{proof}

\paragraph{Complex distortion of the push-forward map.}
Consider a filling-in construction, a point $z$ not on the real
line, and and integer $k$. Let $\Upsilon_{\infty}$ mean the limiting 
branchwise equivalence.
By the push-forward step we mean the following procedure.  
Find the largest $i$ so that $\Upsilon_{j}(z)$ is the same same as
$\Upsilon_{k}(z)$ for all $k>j\geq i$. Assume that $i>1$. Then,
as argued in the proof of Lemma~\ref{lem:807p,1}, $\Upsilon(z)$ is in
the pull-back region. Thus, the push-forward map can applied to find $z_{1}$  
We also put $k_{1}=i-1$. Denote this push-forward map with $\zeta_{1}$.
We can the repeat the push-forward step with $z:=z_{1}$ and
$k:=k_{1}$, and thus construct the sequences up to $z_{l}$ and
$k_{l}$. The construction may end when $i_{l}=0$ or $1$. In the latter
case $z_{l}$ is in the glueing region. 

We can then compose the push-forward maps
$\zeta_{l}\circ\cdots\circ\ zeta_{1}$. An important property
of our construction is that  
\[\Upsilon_{k}=(\hat{\zeta_{l}}\circ\cdots\circ\hat{\zeta_{1}})^{-1}\circ
\Upsilon_{i_{l}}\circ\zeta_{l}\circ\cdots\circ\ zeta_{1}\]
on a neighborhood of $z$. In the future, we will need this lemma:
\begin{lem}\label{lem:q24p,1}
For any $z$, the complex distortion of the corresponding composition
on $z'$ in a neighborhood of $z$
$\zeta_{l}\circ\cdots\circ\ zeta_{1}$ is bounded in a uniform way 
proportionally to the $y$-coordinate of $z_{l}$.
$\Upsilon$.
\end{lem}
\begin{proof}
All maps $\zeta_{i}$ are local extensions of folding branches of
$\Upsilon$. We claim that the complex distortion of the composition at
any point $z'$ in a neighborhood of $z$ where the composition is
defined is bounded proportionally to the sum of the $y$-coordinates of
points $z_{i}$. Indeed, by Lemma~\ref{lem:n6,4}, the complex
distortion of $\zeta_{i}$ can be bounded proportionally to the sum of
$y$-coordinates of $z_{i}$ and $z_{i-1}$. 

But, as in the proof of Lemma~\ref{lem:807p,1}, we argue that the
$y$-coordinates grow exponentially with the push-forward step.  The
lemma follows.    
\end{proof}

\section{Estimates of conformal distortion}  
\subsection{Global description of the construction}
\paragraph{Postulates.}
In the previous section we learned how to realize complex steps of
individual pull-back, chain refinement and filling-in. Thus, we
already know that the construction described in the Branchwise
Equivalences section could be traced by these complex procedures.
We are now ready for our most challenging task of estimating the
conformal distortion of the maps we get. 

We start with an abstract approach by defining an admissible complex
construction.

\subparagraph{An admissible complex construction.}
We begin with a primary complex branchwise equivalence which we will
need in four versions: non-boundary-refined, fully boundary-refined,
and boundary-refined on each side. We still reserve the right to
choose these primary branchwise equivalences suitably fine.  

Then, we proceed to build more branchwise equivalences by these steps
used in an arbitrary order:
\begin{itemize}
\item
A monotone or critical pull-back on a single branch.
\item
A simultaneous chain pull-back. 
\item
A filling-in as described in the previous section.
\end{itemize}

In addition, we assume that the construction is conducted so the
length ratio of any pair of adjacent domains is always uniformly bounded.
Also, we inductively define ``irregularity'' of a point on the real
line as follows. All points of the primary equivalences receive
irregularity $0$. If $\Upsilon_{1}$ is refined by pulling back
$\Upsilon_{2}$, the irregularity at a point is equal to the sum of its
irregularity with respect to $\Upsilon_{1}$ and the irregularity of
its push-forward image relative $\Upsilon_{2}$ increased by $1$ if the
point is an endpoint of a branch being refined and $\Upsilon_{2}$ is
not boundary refined on the side of the push-forward image of the
point.   

We will later add one more assumption, but we need to preparations to
state it clearly. 

We notice that extensions of real branchwise equivalences built in the
Branchwise Equivalences section can be obtained in an admissible
construction. The estimate on the length ratio of adjacent domains
follows directly from Proposition ~\ref{prop:11,2} and the second
property is also provided by the real construction.

From our construction, we see that $D(\zeta_{i})$ are restricted by
two conditions: that they must be inside $D(\xi)$, and also inside the
preimages of $D(\zeta_{i})$.
\paragraph{The tree.}
The complex construction is quite complicated, and trees can be used
to describe and better understand it. We now understand the
construction as a set of branchwise equivalences where each comes with
the prescription for how to build it from other  branchwise
equivalences so that it is possible to ultimately reduce it to the
primary branchwise equivalences. 

A vertex of the tree is the following triple: a point $z$ of the complex
plane, a branchwise equivalence $\Upsilon_{1}$ and its refining map
$G$.  

Each vertex may have up to two daughters: one ``left'' and one ``up''.
If $\Upsilon_{1}(z)$ is not in the pull-back region of $G$, we look at
how $\Upsilon_{1}$ was constructed. If $\Upsilon_{1}$ is primary,
there are no daughters. Otherwise, it is equal to
$G_{0}\circ\Upsilon_{0}$. Then $(z,\Upsilon_{0}, G_{0})$ is the left
daughter, and still there is no up daughter.

If $\Upsilon_{1}(z)$ is in the pull-back region, we look at its
push-forward image $\rho(z)$ in the domain of $\Upsilon_{2}$. We find
the left daughter as in the previous case, and there is an up
daughter, too. If $\Upsilon_{2}$ was obtained as
$G_{3}\circ\Upsilon_{3}$, the up daughter is $(\rho(z),\Upsilon_{3},
G_{3})$. 

Thus, given one vertex a tree can be built according to these rules. 
\subparagraph{Degrees of branches.}

Given a tree, we introduce the  {\em degree} of a folding branch. By 
definition, it is $0$
for all branches of the primary map. Whenever a new central branch is created,
its degree grows by $1$ compared with the degree of the old central branch.
Finally, the degree of a folding branch which is the preimage of some
central branch is equal to the degree of that central branch.  

Clearly, the length of the domain decreases exponentially fast with the degree
and the ratio can be controlled by choosing the primary map. 

\subparagraph{Restrictions on admissible constructions.}
We make two more assumptions about our admissible constructions.
\begin{itemize}
\item 
Given any $\Upsilon$ and its refining map $G$, the tree built from
$(z,\Upsilon,G)$ is finite except for $z$ from a set of zero measure. 
\item
In any tree of the construction if we follow a branch up, the degrees of
folding branches being refined form a non-increasing sequence. The degrees of
consecutive branches can be equal only if  they belong to consecutive
branchwise equivalences in a filling-in step.  
\end{itemize}

From now on when we speak of admissible constructions, we mean that
these two properties hold as well.

\subparagraph{Push-forward along vertical branches.}
In the previous section, we defined the push-forward map on any
component of the pull-back region. Now, the push-forward map may be
associated with a vertical edge in the tree of some point. Indeed,
choose the triple $(Upsilon, G, z)$ with any $z$ in the domain of the
push-forward map, and the map itself is defined. Analogously, the
definition can be extended so that we can push-forward along any
vertical branch. The definition, which now depends on the choice of
some tree, goes simply by composing the push-forward maps which
correspond to consecutive vertical edges. The map is defined wherever
the composition is.
  
\paragraph{Small diamonds.}
\subparagraph{Bounded sizes of close neighborhoods.}
\begin{lem}\label{lem:n7,1}
The sizes of close neighborhoods in an admissible construction are
uniformly bounded away from $0$.
\end{lem}
\begin{proof}
This is clearly true of primary branchwise equivalences whose close
neighborhoods have unit size. If we examine close neighborhoods of
branches obtained in a pull-back operation, we notice that they
certainly contain preimages of the close neighborhoods of the branches
being pulled back intersected with the large diamond around the branch
being refined. From this, it follows by induction that for any branch
with stopping time $s$ all points whose forward orbits stay in large
diamonds around intermediate images are in the close neighborhood.
This set contains a diamond of fixed size by Proposition~\ref{prop:11,1}.  
\end{proof}

These diamonds of fixed size will be called {\em small diamonds.}

We notice that, in addition to being contained in close neighborhoods,
small diamonds have this property:

If $z$ is in the small diamond of a branch created by refinement of a
branchwise equivalence $\Upsilon$ by a refining map $G$,  any
push-forward map constructed for the vertical branch starting at 
$(\Upsilon, G, z)$ is defined on the whole small diamond.  

Actually, the defining formula of close neighborhoods may be regarded
as a corollary from this ``push-forward'' property when the push-forward
goes all the way up to the primary equivalence. 

\subparagraph{Bounded conformal distortion of marking.}
When we defined complex marking, the question was left unanswered of
the complex distortion of the modified map $A$. Now, in view of 
Fact~\ref{fa:n7,1} and Lemma~\ref{lem:n7,1} we see that the
modification  can be done in a uniformly quasiconformal fashion.   

\subparagraph{The choice of primary branchwise equivalences.}
We can choose 
the primary branches short enough so that the following is satisfied:
\begin{itemize}
\item
Any complexified branchwise equivalence coincides with the primary branchwise
equivalence except on inside a diamond neighborhood of the unit interval. 
The size of that neighborhood can be made arbitrarily small by choosing 
suitably short branches in the primary map.
\item
Suppose that a branch $\zeta$ is being refined. Subsequently, consider 
a chain refinement which involves any branch  created inside $\zeta$, possibly
after many steps. Then, the corresponding refining map is the identity outside
the small diamond neighborhood of $\zeta$, unless a branch inside $\zeta$ 
adjacent to an endpoint of $\zeta$ is involved in the chain refinement.

In particular, if the original refinement was a pull-back of a 
boundary-refined
map, no future refinements inside of $\zeta$ will ever affect the region 
outside of the diamond.   
\end{itemize}

\subparagraph{Quasiconformal estimates for the refining map.}

Clearly, the quasiconformal distortion of the map is null on the
trivial 
region,
while on the pull-back region it strictly depends on the properties of
the maps
being pulled back. We want the following estimate:

The quasiconformal distortion is uniformly bounded on the glueing region.

Let us think for a moment
that only monotone branches are being refined. Then, it is sufficient to choose
a suitably fine primary branchwise equivalence. Indeed, we showed in section
3 that the glueing operations only contribute a uniformly bounded distortion.
The only issue is to show the potentially unbounded distortion coming from 
the pullback itself is supported inside the pull-back region. We noted the 
the potentially unbounded distortion of the map being pulled back is supported
inside a diamond, which can be made tiny by choosing the primary equivalence
appropriately. Thus, it will also be a tiny diamond after the monotone
pull-back, and we can choose the constants so that in fact it fits inside the
pull-back region. The critical pull-back poses a problem, though. There will
be a part of the preimage of the diamond which sticks out.

This is why we construct marked maps in a special way, so that we do
not refine the primary 
branches whose preimages are going to be imaginary. Then, the same
argument which we have used for monotone pull-backs still applies.  

\subsection{Estimates}
\paragraph{Rough distortion.}
We will estimate quasiconformal distortion in terms of a ``combinatorial''
object that we call {\em rough distortion}. The definition is as follows:
\begin{defi} \label{defi:n1}
For any branchwise equivalence, we define an integer valued function on the
plane. For the primary branchwise equivalence it is identically $0$. For a 
refining map, it is $0$ in the trivial region, the same as at the corresponding
points of the pulled-back map in the pull-back region, and $1$ in the glueing 
region. Finally, after a pull-back operation, the value of the function is the
sum of its value for the map being refined and for the refining map at the 
image.

The function so defined is called the {\em rough distortion}.
\end{defi}

\subparagraph{Remark:} Thus, precisely speaking, the rough distortion depends
not just on the branchwise equivalence, but also on the way it was obtained
(although this way is in fact unique in our construction.) Since all our 
arguments are recursive, that makes no difference.  

\paragraph{Complex distortion bounded by rough distortion.}
Here is an important lemma.

\begin{prop} \label{prop:n1,2}
There is a function $Q(n)$ so that for any branchwise equivalence if 
the rough distortion at a point is $n$, the quasiconformal distortion is
bounded by $Q(n)$.
\end{prop}
\begin{proof}
We will prove that by induction with respect to the rough distortion.
Fix your attention on the map being refined and some point $z$. We look for 
the last refinement step that changed the map in a neighborhood of $z$.
For that step, $z$ cannot be in the trivial region. In $z$ is in the glueing 
region, we are done. Indeed, the rough distortion must have grown by $1$ 
compared with the map being refined. So, this cannot happen at all at the 
initial step of the induction (rough distortion equal to $0$), otherwise
an estimate follows from the fact that the complex distortion of the glueing 
map in the glueing region is uniformly bounded. 

So, the real problem
occurs if $z$ is in the pull-back region. Then, we consider its push-forward
image $z'$ in some branchwise equivalence. If $z'$ again is in the
pull-back region of some refinement, we can iterate the procedure.
Thus, we get a sequence of points $z_{0}:=z,\ldots,z_{k}$ of images by
consecutive push-forward maps, and $\Upsilon_{k}(z_{k})$ is no longer
in the pull-back
region. Let $\zeta_{0},\ldots,\zeta_{k-1}$ denote the consecutive
push-forward maps. Since the glueing regions are open, the composition  
$\zeta_{k-1}\circ\cdots\circ\zeta_{0}$ is defined on an neighborhood
of $z$.
To complete the proof of the proposition, it will be enough if show
that the complex norm of this composition is uniformly bounded.    
Indeed, by the properties of the pull-back region, the branchwise
equivalence on a neighborhood of $z$ is given by
\[(\hat{\zeta_{k-1}}\circ\cdots\circ\hat{\zeta_{0}})^{-1}\circ
\Upsilon_{k}\circ\zeta_{k-1}\circ\cdots\circ\zeta_{0}\; .\] 
As we already noted, the complex distortion of $\Upsilon_{k}$ at
$z_{k}$ can be bounded from induction. 

Every point $z_{i}$ for $i<k$ can be associated with a branch, namely the only
in whose small diamond it is. We call this branch $\xi$.

Next, we seek out sequences of critical pull-backs which correspond to
a filling-in operation. We regard the push-forward map which
corresponds to the filling-in (see previous section) as just one map.  

So, we get subsequences indexed by $i_{j}$. Maps $\zeta_{i_{j}}$ are
now of three types: local extensions of monotone and folding branches,
and push-forward maps of the filling-in.

Then, we observe that the complex distortion of the whole composition at $z$  
is bounded in a uniform fashion proportionally to the sum of lengths
of $\xi_{i_{j}}$. Indeed, look at the $\zeta_{i_{j}}$ at $z_{i_{j}}$. 
The height of $z_{i_{j}}$ 
relative the domain $\xi_{i_{j}}$ as well of $z_{i_{j+1}}$ relative the domain
$\xi_{i_{j+1}}$ is bounded for any $i_{j+1}<k$ by virtue of both points being
in their respective small diamonds. 

The map $\zeta_{i_{j}}$ is nothing else but the local extension of a
branch being refined. If it is monotone, or of the
form
\[ h_{1}\circ Q\circ h_{2}\]
where $h_{1}$ and $h_{2}$ are monotone so that Lemma~\ref{lem:n6,4}
can be used to bound the complex distortion of
the extensions of monotone branches. The result is that the distortion
of $\zeta_{i_{j}}$ is bounded is bounded by the sum of lengths of the
domains of $\xi_{i_{j}}$ and $\xi_{i_{j+1}}$. If $\zeta_{i_{j}}$ a
filling-in push-forward, the estimate follows directly form
Lemma~\ref{lem:q24p,1}.      

This reasoning does not apply to $\zeta_{k-1}$, but its contribution
is also uniformly bounded. Thus, we only need to sum up the
contributions for all $i$ to the bound proportional to the sum of
lengths of the domains.

This reduces the problem to the real line.

Next, we pick a subsequence of $j$, which we denote with $j_{l}$. 
An index $j$ enters this subsequence unless
$\zeta_{i_{j}}$ is monotone. We claim the sum of lengths of all
domains is bounded proportionally to the sum of lengths of domains of
$\xi_{i_{j_{l}}}$. Indeed, consider the domains of $\xi_{i_{m}}$ with
$j_{l}<m\leq j_{l+1}$. Since $\zeta_{i_{m}}$ are monotone except for
$m=j_{l+1}$, the lengths of $\xi_{i_{m}}$ increase exponentially with
$m$. Thus, the total is bounded by the last term, which gives our
claim.

Finally, the lengths of the domains of $\xi_{i_{j_{l}}}$ grow
exponentially with $l$ as a direct consequence of admissibility of the
construction, namely, that the degrees of folding branches decrease up
any vertical branch the tree unless consecutive vertices belong to the
same pull-back operation. 

The proposition follows.
\end{proof}

Thus, it remains to show that the rough distortion is uniformly bounded, which
is what we do next.

\paragraph{Boundedness of the rough distortion.}

\begin{lem} \label{lem:n1,1}
Inside the small diamond neighborhoods the rough distortion is $0$.
\end{lem}
\begin{proof}
By the definition, the rough distortion at a point inside the large diamond 
is the same as the rough distortion at its image by push-forward map. But the
small diamond was defined by the property that the push-forward map can be 
iterated all the way, and for the primary branchwise equivalence the rough 
distortion is $0$.
\end{proof}

We look for the simplest branchwise equivalence so that the rough distortion 
at a point $z$ is $k$. That means that when the branchwise equivalence was
created, neither the pulled-back map nor the refined map had points with 
rough distortion $k$. Consider the map being refined. Clearly, the image of 
$z$ is in the glueing region. Then, look at the branch directly below
$z$. Call this branch $\zeta$.
If $z$ is above the boundary of two branches, take any of them.
Observe that $z$ cannot be above the Cantor set, since such points
are fixed by subsequent construction.

Since the image of $z$ was in the glueing region, the height of $z$
with respect to $\zeta$ is bounded away from $0$ and infinity. Then look for 
the refinement step when $\zeta$ was created. Unless $\zeta$ was adjacent to
the endpoint of the branch then being refined, $z$ was in the small diamond.
If $z$ is in the small diamond of the branch being refined, we can push both 
$z$ and $\zeta$ forward and look at the corresponding
objects (we continue to call them $z$ and $\zeta$.) Then, we can repeat the 
procedure. Thus, we arrive at one of two possible outcomes: either we can push
forward to the primary map, or at some point $\zeta$ is adjacent to the 
boundary of the branch being refined, and moreover the $\zeta$ is 
outside of the small diamond neighborhood of that branch. In the first 
situation the rough distortion at $\zeta$ is $0$, so $k$ is $1$ and we are 
done. Let us consider the branch being refined and call it $\zeta_{1}$.
 Clearly, $Z$ must also be 
the end of the whole chain, because otherwise the adjacent refinements would
both be boundary-refined and there would be branch adjacent to $z$. 
The point $z$ must be very close to an endpoint od $\zeta_{1}$ which we call
$Z$ compared with the length of $\zeta_{1}$.
How close again depends on the 
choice of the primary map. The rough distortion at $z$ is now at least $k-2$.
We use the same argument to $z$ and $\zeta_{1}$. That means that we either can
push them forward to the primary map, or a push-forward image coincides with
an endpoint of the branch being refined, in which case the rough distortion
distortion may drop by $1$. However, our construction ensures that among
the push-forward images of any point at most two are endpoints of chains.

Thus, we can get at most one repetition of this situation. So, $k\leq 3$.

By Proposition~\ref{prop:n1,2} this concludes the proof of the main theorem.  

\subsection{Invariance of the Collet-Eckmann condition}\label{sec:5,3}

The construction of \cite{jak1} provides maps with
absolutely continuous invariant measures which are all basic. 
Such maps constitute a  positive  measure
set of parameter values for typical one-parameter
families of S-unimodal maps. As proved
by Benedicks and Carleson (see \cite{beca})  the same is true
for maps satisfying  Collet-Eckmann condition which can be written as
\[ Df^{n} (0)  > ab^{n} a>0, b > 1 \]
for every $n>0$. 

Theorem1 implies that: 

\paragraph{Corollary}
\begin{em}
For maps from $\cal C$ with basic dynamics, the  Collet-Eckmann
condition is a topological invariant. 
\end{em}

\subparagraph{Proof of the Corollary}
The basic construction of \cite{yours} results in a partition of
$(-q,q)$ into domains of monotone branches $f_{i}$ which
are uniformly extendable and so are all their compositions.

Next, we notice that if $f_{i}=f^{n_{i}}$, then all compositions 
$f^{j},\:\: j\leq n_{i}$ are also extendable from the domain of
$f_{i}$. Indeed, the ``space'' around the image of $f^{j}$ is the
preimage of the space around the image of $f^{n_{i}}$ by a negative
Schwarzian map, hence it is large.  

Thus, the derivatives in the Collet-Eckmann condition are approximated
up to a multiplicative constant by ratios of lengths of dynamically defined
intervals.

But the qs conjugacy is also H\"{o}lder continuous and so exponential 
decreasing of such ratios is preserved.

\section{Renormalizable polynomials}\label{sec:6}
In this section, we extend our results to a certain class of renormalizable
S-unimodal maps, including some infinitely renormalizable ones for which the 
result is new even in the polynomial case.

\subsection{Statement of the problem}
\paragraph{Restrictive induced maps.}
\begin{defi}\label{defi:n5,3}
Suppose that $\phi$ is a preferred induced map or is unimodal, that
is, consists of one folding branch.
 Suppose that the critical value of
$\phi$ is in the domain of a folding branch $\psi$, moreover, under the 
iterations of $\phi$ the critical orbit forever stays inside the domain of 
$\psi$. If that happens, we say that  $\phi$ is a restrictive induced map. 
\end{defi}

\begin{lem}\label{lem:n5,1}
If $\phi$ is a restrictive induced map, then the underlying $f$ has a 
restrictive interval. If $h_{1}$ 
means the natural diffeomorphism from the domain of $\psi$ onto the central 
domain of $\phi$, then $h_{1}\circ\phi$ is the first return map on this 
interval. 
\end{lem}
\begin{proof}
Standard.
\end{proof}

\subparagraph{Remark.} Suitable induced maps mentioned in the
introduction are restrictive induced maps in the sense of
Definition~\ref{defi:n5,3}. 

\paragraph{Renormalization.}
Let $\phi$ be a restrictive induced map, $I$ be the restrictive interval from
Lemma ~\ref{lem:n5,1}, and $f_{1}$ be the first return map onto $I$. Which 
rescale $I$ affinely so that it becomes the unit interval. The first return 
map gives us some unimodal endomorphism of $[0,1]$, which we will also call 
$f_{1}$. The basis of argument is this:
\begin{fact}\label{fa:n5,1}
Under our assumptions, $f_{1}\in {\cal C}$. Moreover, if 
$f_{1}=h_{1}(x^{2})$ the distortion of $h_{1}$ is bounded in a uniform way.
\end{fact}

This follows directly from Theorem 1 in \cite{miszczu}. An
equivalent result of \cite{blyub} should also be noted.

\subparagraph{Matching sequences of restrictive induced maps.}
Let $\varphi_{0}, \varphi_{1}, \varphi_{\omega}$ be a sequence of
induced maps, either finite or infinite. All of them are restrictive
except for the last one, $\varphi_{\omega}$, if it exists. Also, we
assume that $\varphi_{0}$ is an induced map on the standard domain
$[-q,q]$ for some interval map $f$. 

We say that a sequence which satisfies all these properties is a {\em
matching sequence of restrictive induced maps} if $\varphi_{i+1}$ and 
$\varphi_{i}$ are related as follows:

associate a map $f_{1}$ with $\varphi_{i}$ as in the previous paragraph.
Then $\varphi_{i+1}$ is an induced map on the standard domain for $f_{1}$.

\paragraph{The main result.}
\begin{prop}\label{prop:n5,1}
Let $\varphi_{0},\ldots,\varphi_{\omega}$ be a matching sequence of restrictive
induced maps for $f$, and $\hat{\varphi_{0}},\ldots,\hat{\varphi_{\omega}}$ be
an analogous sequence for a topologically conjugate map $\hf$. We
allow $\omega$ to be infinite. 
 
We assume that all $\varphi_{i}$ with the possible exception of
$\varphi_{\omega}$ are regular. 
Suppose, 
further, that admissible complexified branchwise equivalences
$\Upsilon_{i}$ are given 
between $\varphi_{i}$ and $\hat{\varphi}_{i}$ for which
small diamonds can be chosen with uniform size, and which are
uniformly quasiconformal. 

If $\omega$ is finite, then $\Upsilon_{\omega}$ is assumed to be
conjugacy.

Then, $f$ and $\hf$ are quasisymmetrically conjugate. Moreover, the qs
norm of the conjugacy is bounded by a continuous function of
the small diamond size and the supremum of qc norms. 
\end{prop}

The rest of this section will devoted to the proof of Proposition
~\ref{prop:n5,1}. Before we tackle the proof, we would like to give a
couple of simple corollaries.
\subparagraph{Corollaries.}
First of all, Proposition ~\ref{prop:n5,1} can be used in the
``basic-renormalizable'' case. This case can be characterized by the 
requirement that all maps of a matching sequence of restrictive
induced maps can be obtained by the basic construction, that is, in
the process of their construction the critical value of the
intermediate preferred induced maps never falls into a folding domain. 
By the results of previous sections, complexified branchwise
equivalences can be constructed which satisfy the assumptions of
Proposition ~\ref{prop:n5,1}. So, our main theorem can be extended on
the basic-renormalizable case.

This also includes a ``basic-finitely renormalizable'' case in which
the non-renormalizable map obtained in the last box is also assumed to
be basic. Then, we use Theorem 1 and Proposition~\ref{prop:n5,1} with
finite $\omega$ to prove quasisymmetric conjugacy.

Also, the case of Feigenbaum, or ``bounded type'', maps considered in 
\cite{miszczu} is  reduced to Fact ~\ref{fa:n5,1}. The
``bounded type'' assumption means that the number of branches of all
maps from the matching sequence of restrictive induced maps is
uniformly bounded. But then, their sizes must be uniformly comparable
as an easy consequence of Fact ~\ref{fa:n5,1}, and maps $\Upsilon_{i}$
can be constructed ``by hand'' to satisfy the hypotheses of
Proposition ~\ref{prop:n5,1}. 
 
\subsection{A single matching step}
Not surprisingly, the idea of the proof of Proposition~\ref{prop:n5,1}
is to somehow imprint the structure given by $\Upsilon_{i+1}$ into the
restrictive which lives somewhere in $\Upsilon_{i}$ and continue with
this process. To do this in a way that, given our results about
admissible constructions, will automatically ensure a bounded qc norm
of the result, we need to ``prepare'' $\Upsilon_{i}$ for this
operation. Our matching step proposition is about that. 

\paragraph{The matching step proposition.}
\begin{prop}\label{prop:n5,2}
Suppose that we have a regular restrictive induced map $\varphi$
an the corresponding admissible complexified branchwise equivalence 
$\Upsilon$ with uniformly large small diamonds. 

Then, an admissible complex boundary-refined branchwise equivalence 
$\Upsilon'$ can be built which satisfies the following conditions:
\begin{itemize}
\item
The quasiconformal norm of $\Upsilon'$ is bounded by a uniform
function of the qc norm of $\Upsilon$.
\item
All branches are monotone and map onto the fundamental inducing domain
of the first return map on the restrictive interval.
\item
The marked set comprises the complement of the union of the domains of
the branches.
\end{itemize}
\end{prop}

\subparagraph{Derivation of Proposition~\ref{prop:n5,1}.}
We will show how Proposition ~\ref{prop:n5,1} follows from
Proposition ~\ref{prop:n5,2}. 

If the matching sequence is infinite, we choose $\omega$ in an
arbitrary fashion. We consider the following complex construction.

Maps $\Upsilon'_{0},\ldots,\Upsilon'_{\omega-1}$ and a
boundary-refined version of $\Upsilon_{\omega}$, called
$\Upsilon'_{\omega}$,  are primary where
$\Upsilon'_{i}$, $i<\omega$, means the map which corresponds to 
$\Upsilon_{i}$ by Proposition~\ref{prop:n5,2}. 

We construct a sequence $\Upsilon^{i}$ defined inductively.
$\Upsilon^{\omega}$ is equal to $\Upsilon'_{\omega}$. For $0\leq
i<\omega$, $\Upsilon^{i}$ is form by the pull-back of $\Upsilon^{i+1}$
onto all branches of $\Upsilon'_{i}$. It should be noted that this
operation can be realized as a simultaneous chain monotone refinement,
usually with many chains.

Since this is an admissible construction, the qc norm of
$\Upsilon^{0}$ is bounded uniformly as a function of the maximum of norms of  
$\Upsilon'_{i}$, hence of $\Upsilon_{i}$. 

If $\Upsilon_{\omega}$ was a
conjugacy, $\Upsilon^{0}$ is, too. Otherwise, the matching sequence is
infinite. We notice that the $\Upsilon^{0}$ in that case coincides
with the conjugacy except on the interior of all preimages of the
restrictive interval of $\varphi_{\omega}$. But $\omega$ can be chosen
arbitrarily, and, as it grows, the complement of this set grows to a
dense set (tends to the whole interval in the Hausdorff distance
uniformly with $\omega$). Hence, the corresponding maps $\Upsilon^{0}$
tend to the conjugacy on the line in the $C^{0}$ topology. Even though
it is not obvious that they converge everywhere, they are a normal
family, since they are uniformly continuous and all identical except
on a compact set. 

Proposition ~\ref{prop:n5,1} follows.

\paragraph{An outline of the construction.}
Let $\psi$ mean the folding branch which fixes an image of the
restrictive interval.

We will first describe how to construct $\Upsilon'$ on
the real line. We will do so in familiar terms of pull-backs so that
the complexification of this procedure will be easy.

Suppose that $\varphi$ is not unimodal.

First, we want to pull the structure defined by $\Upsilon$ into
the domain of the branch $\psi$. We notice that each point
of the line which is outside of the restrictive interval will be
mapped outside of the domain of $\psi$ under some number of
iterates of $\psi$. We can consider sets of points for which the 
number of iterates required to escape from the domain of $\psi$ is
fixed. Each such set clearly consists of two intervals symmetric with
respect to the critical point. The endpoints of these sets form two
symmetric sequences accumulating at the endpoints of the restrictive
interval, which will be called {\em outer
staircases}. Consequently, the connected components of these sets will
be called steps. 
 
This allows us to construct an induced map
from the complement of the restrictive interval in the domain of
$\psi$ to the outside of the domain $\psi$ with branches
defined on the steps of the outer staircases.  That means, we can
pull-back $\Upsilon$ to the inside of the domain of $\psi$. 

Next, we construct the {\em inner staircases}. We notice that every
point inside the restrictive interval but outside of the fundamental
inducing domain inside it is mapped into the fundamental inducing
domain eventually. Again, we can consider the sets on which the time
required to get to the fundamental inducing domain is fixed, and so we
get the steps of a pair of symmetric inner staircases.

So far, we have obtained a branchwise equivalence which has one
indifferent domain equal the restrictive interval and besides has
extendable monotone and folding branches. Denote it with
$\Upsilon^{1}$.  The folding branches are all
preimages of $\psi$. We now refine the folding branches. 

This can be done in the usual filling-in way.

We conclude with refinement of remaining monotone branches analogous
to the final refinement step in the basic case.
Since there are no folding branches left, we can 
destroy all monotone branches. Thus, we will be left with indifferent
branches only, all of which are certain preimages of the restrictive 
interval. So, at least topologically, we obtain a good candidate for
$\Upsilon'$ on the line. 
 
In the case when $\varphi$ is unimodal, the outer staircase cannot
be constructed. Instead we build the inner
staircases twice. The first step is as described. For the second step,
we notice that the restrictive interval is the same as the fundamental
inducing domain of $\psi$. So, the inner staircases can be built
again.

This completes the real description of the matching step. What remains
is to define the complex version of this procedure and do estimates.

\paragraph{Outer staircases.}
Unless we explicitly indicate otherwise, the assumption is that
$\varphi$ is not unimodal.

\subparagraph{Pulling-in outer staircases from far away.}
Suppose that the domain of $\psi$
is very short compared with the length of the the domain of
$\varphi$. This means that the domain of $\psi$ is extremely
large compared with the restrictive interval. This unbounded situation
leads to certain difficulties and is dealt with in our next lemma.

\begin{lem}\label{lem: n5,1}
One can construct a map $\Upsilon^{1}$ which is an admissible complex
branchwise equivalence and its
qc norm, as well as the sizes of its small diamonds are uniformly
related to the analogous estimates for $\Upsilon$. In addition, an
integer $i$ can be chosen so that the following conditions are satisfied:
\begin{itemize} 
\item
The functional equation
\[ \Upsilon\circ\psi^{j}=\hat{\psi}^{j}\Upsilon^{1} \; \]
holds for any $0\leq j\leq i$ whenever the left-hand side is defined.
\item
The length of the interval which consists of points whose $i$
consecutive images by $\psi$ remain in the domain of $\psi$ forms a
uniformly bounded ratio with the length of the restrictive interval. 
\end{itemize}   
\end{lem}
\begin{proof}
We rescale affinely so that the restrictive intervals become $[-1,1]$
in both maps. Denote the domains of $\psi$ and $\hat{\psi}$ with $P$
and $\hat{P}$ respectively.
Then, $\psi$ can be represented as $h(x^{2})$ where
$h''/h'$ is very small provided that $|P|$ is large. We can assume that $|P|$
is large, since otherwise we can take $\Upsilon^{1}:=\Upsilon$ to
satisfy the claim of our lemma.
Thus, assuming that $|P|$ is large enough, we can find a uniform $r$
so that the preimages
of  $B(0,r)$ by $\psi$, $\hat(\psi)$ and $z\rightarrow z^{2}$
are all inside $B(0,r/2)$. Also, we can have $B(0,r)$ contained in the small 
diamond around the domain of $\psi$. Next, we choose the largest $i$ so that
$[-r,r]\subset\psi^{-i}(P)$ . 

Then, we change $\psi$ and $\hat{\psi}$. We will only
describe what is done to $\psi$. Outside of $B(0,r)$, $\psi$
coincides with its standard extension. Inside
the preimage of $B(0,r)$ by $z\rightarrow z^{2}$ it is $z\rightarrow
z^{2}$. In between, it
can be interpolated by a bounded distortion smooth 2-1 local
diffeomorphism. We leave to the reader to convince himself it is
possible to construct such a map. Also, look up \cite{miszczu} where a
similar situation is considered. The modified extension will be
denoted with $\psi'$.

Next, we pull-back $\Upsilon$ by $\psi'$ and $\hat{\psi}'$ exactly $i$
times. That is, if $\Upsilon_{0}$ is taken equal to $\Upsilon$, then
$\Upsilon_{j+1}$ is $\Upsilon$ refined by pulling-back $\Upsilon_{j}$
onto the domain of $\psi$. This perhaps requires a little further
clarification, since in Section~\ref{sec:3} we only defined the
pull-back by branches and under the assumption that the critical
branch was in a monotone domain. Here, we mean the the simple pull-back is
obtained by the same formula that would be used to pull-back by
$\psi$, only $\psi$ is replaced by $\psi'$. Note that no extra
marking is required as the critical value of $\psi$ is at $0$ and $0$
is preserved automatically. This determines the map inside the small diamond. 
Outside of it, the refining map is corrected in the usual way. Note
that the middle branch of the map so constructed is not ``folding''
since it is 
of degree $2^{i}$ rather than quadratic. Hence, it does not satisfy our
definition of the folding branch and must be considered indifferent.
$\Upsilon_{i}$ constructed in this way can be taken as $\Upsilon^{1}$.

Now we need to check whether $\Upsilon^{1}$ has all the properties
claimed in the Lemma. To see admissibility and the functional equation
condition, we note that all branches of any $\Upsilon_{j}$ and their
small diamonds are in the region where $\psi$ coincides with
$\psi'$. Thus, the same arguments as in Section~\ref{sec:3} can be
used to prove admissibility, and the functional equation is also
evidently true. 

The last condition easily follows  from the fact that $r$ can be chosen in a
uniform fashion. 

So, what remains is estimates of the qc norm $\Upsilon_{i}$. 
First, we note that the qc distortion of $\Upsilon_{j}$ for any $j\leq
i$ at points not inside $B(0,r)$ is bounded as a uniform function of
the qc norm of $\Upsilon$. We notice that $\Upsilon_{j}$ in the
complement of $B(0,r)$ is the pull-back of $\Upsilon$ by unmodified
$\psi$. So, as usual, we can use the fact that the distances of
push-forward images from the line grow exponentially, thus the total 
distortion added is bounded. 

Points inside $B(0,r)$ are pull-backs of points outside of $B(0,r)$ by
$\psi'$. But, $\psi'$ is conformal inside
$B(0,r)\setminus\psi^{,-1}(B(0,r))$ and quasiconformal inside
$B(0,r)$. Also, only one push-forward image is inside the region where
the map is not conformal. So, again, only bounded distortion is
acquired.
\end{proof}

\subparagraph{Comment.}
One should be aware that the situation handled in Lemma
~\ref{lem:n5,1} is not a bounded pull-back situation. The proportions
of the preimages of the domains on the last step constructed may be
arbitrarily different from the proportions on the zeroth step. For
example, even if $\Upsilon$ is not boundary-refined, it is not
true that the preimage of the outermost domain constitutes any fixed
part of the last step. 

\paragraph{The staircase construction.}
We take $\Upsilon^{1}$ obtained in Lemma~\ref{lem:n5,1} and restrict
our attention to its restriction to the real line, denoted with $\upsilon_{1}$.
We rely on the fact that $\upsilon_{1}$ is a quasisymmetric map and
its qs norm is uniformly bounded in terms of the quasiconformal norm
of $\Upsilon^{1}$. 

For a while, we will be working with real methods.

\subparagraph{Completion of outer staircases.}
We will construct a
real map $\upsilon_{2}$ from the domain of $\varphi$ to the domain
of $\hat{\varphi}$ with following properties:
\begin{itemize}
\item
The map $\upsilon_{2}$ coincides with $\upsilon_{1}$ outside of the domain
of $\psi$. Also, it satisfies 
\[ \upsilon_{2}\circ \psi^{j} = \hat{\psi}^{j}\circ\upsilon_{2} \]
on the complement of the restrictive interval provided that
$\psi^{j}$ is defined.
\item
Inside the restrictive interval, it is the ``inner staircase
equivalence'', that is, all endpoints of the inner staircase steps are
mapped onto the corresponding points.
\item
Its qs norm is uniformly bounded as a function of the qc norm of 
$\Upsilon$.
\end{itemize}

Outer staircases constructed in Lemma\ref{lem:n5,1}
connect the boundary points of the domain of $\psi$  to
the $i$-th steps which are in the close neighborhood of the
restrictive interval. Also, the $i$-th steps are the corresponding 
fundamental domains for the inverses of $\psi$ in the proximity of the
boundary of the restrictive interval. 

From Fact~\ref{fa:n5,1},  the derivative of $\psi$ at the boundary of
the restrictive interval is uniformely bounded away from one.

Then, it is straightforward to see that the equivariant correspondence
between infinite outer staircases which uniquely extends
$\upsilon_{1}$ from the $i$-th steps is uniformely qs.

Inside the restrictive interval, the map is already determined on the
endpoints of steps, and can be extended in an equivariant way onto each
step of the inner staircase.

\subparagraph{Re-complexification.}
We want to construct an admissible complex extension of $\upsilon_{2}$ which
is regarded as a branchwise equivalence. This can readily be done by 
Lemma~\ref{lem:n6,5}. 
The result will be called $\Upsilon^{2}$. Note that $\Upsilon^{2}$ is
``primary'' in the sense that all settling times are $1$. As to the
stopping times, they have been defined by the refinement procedure
outside of the restrictive interval. However, we also want to regard
steps of the inner staircase as domains of branches. The stopping time
on a step is going to correspond to the iterate of $\psi$ which maps
this step onto the fundamental inducing domain inside the restrictive 
interval.  

\subparagraph{Remarks on the staircase construction.}

The map $\Upsilon^{2}$ represents the first important step of matching
in the case when $\varphi$ is not unimodal. 
We have built both inner and outer staircases and they fit together.
Moreover, we introduced the structure of an induced map in the outer
staircase, i.e. it is divided into the domains of monotone branches
and folding branches which are copies of $\psi$. In
future, they will be refined and eventually all taken out.

\paragraph{The construction of $\Upsilon'$ when $\varphi$ is unimodal.}

In this case, the construction is quite elementary. There is no outer
staircase, so only the inner staircase is considered. The real map
$\Upsilon_{2}$ maps the steps of one inner staircase onto the
corresponding steps. Then, it is extended beyond the restrictive
interval so that it is affine outside of an interval twice its size,
and is uniformly quasisymmetric. Then, regard it is a boundary-refined
branchwise equivalence, and construct $\Upsilon^{2}$ as its admissible
extension by Lemma~\ref{lem:n6,5}.

Next, we construct the map $\Upsilon_{3}$ which is completely
analogous to $\Upsilon_{2}$, except that it now map the inner
staircase inside formed by preimages of the fundamental inducing
domain of $\varphi$ and not of the first return map on its restrictive
interval. Note that the branches of $\Upsilon_{3}$ map in a monotone
fashion onto the restrictive interval. Thus, we can pull-back
$\Upsilon_{2}$ on them, which a usual chain monotone pull-back. The
result is $\Upsilon'$. That it has the desired properties is clear. 

\paragraph{Final filling of the outer staircase.}
It remains to construct $\Upsilon'$  in the
case when $\varphi$ is not unimodal. 

We need to fill all monotone and secondary folding
branches left outside of the restrictive interval.

\subparagraph{Filling-in of secondary folding branches.}
We perform a filling-in of the folding branches of $\Upsilon_{3}$. 
This only requires one-time marking of the primary branchwise equivalence,
hence presents no problem.

\subparagraph{Final refinement.}  
Then, we apply the final refinement construction to fill all monotone 
branches. Since there is no obstacle presented by
folding branches, the final refinement can be continued until all
monotone branches have disappeared in the limit.

So, we have obtained $\Upsilon'$ in the non-unimodal case.
Its topological properties are evident, and the fact that the qc norm
is suitably bounded follows from out previous estimates.

Proposition~\ref{prop:n5,2} has been proven.

\end{document}